\documentclass[12pt]{amsart}
\usepackage{amssymb}

\usepackage{amscd}

\newcommand{\BR}{{\mathbb{R}}}
\newcommand{\BC}{{\mathbb{C}}}

\newcommand{\BQ}{{\mathbb{Q}}}

\newcommand{\BG}{{\mathbb{G}}}

\newcommand{\BH}{{\mathbb{H}}}
\newcommand{\gD}{\Delta}

\newcommand{\gC}{\Gamma}

\newcommand{\gS}{\Sigma}

\newcommand{\ga}{\alpha}

\newcommand{\SL}{\text{SL}}
\newcommand{\GL}{\text{GL}}
\newcommand{\PSL}{\text{PSL}}

\newtheorem{prop}{Proposition}[section]
\newtheorem{thm}[prop]{Theorem}
\newtheorem{lem}[prop]{Lemma}
\newtheorem{cor}[prop]{Corollary}

\theoremstyle{definition}
\newtheorem{Ack}[prop]{Acknowledgments}

\newtheorem{defn}[prop]{Definition}
\newtheorem{rem}[prop]{Remark}
\newtheorem{exam}[prop]{Example}

\begin{document}
\author{E. Breuillard, T. Gelander}
\date{\today}
\title{A topological Tits alternative}
\maketitle

\begin{abstract}
Let $k$ be a local field, and $\Gamma \leq \text{GL}_{n}(k)$ a linear group
over $k$. We prove that either $\Gamma $ contains a relatively open solvable
subgroup, or it contains a relatively dense free subgroup. This result has applications 
in dynamics, Riemannian foliations and profinite groups. 
\end{abstract}


\section{Introduction}

In his celebrated 1972 paper \cite{Tits0} J. Tits proved the following
fundamental dichotomy for linear groups:
{\it Any finitely generated\footnote{In characteristic zero, one may drop the assumption that the group is finitely generated} linear group
contains either a solvable subgroup of finite index or a non-commutative free subgroup.}
This result, known today as ``the Tits alternative'', answered a conjecture of Bass and Serre and was an 
important step towards the understanding of linear groups. The purpose of the present paper
is to give a topological analog of this dichotomy and to provide various
applications of it. Before stating our main result, let us reformulate Tits' alternative
in a slightly stronger manner. Note that any linear group $\gC\leq\GL_n(K)$ has a Zariski
topology, which is, by definition, the topology induced on $\gC$ from the Zariski topology on $\GL_n(K)$.

\begin{thm}[Tits' alternative]\label{Tits-alternative}
Let $K$ be a field and $\Gamma $ a finitely generated subgroup of GL$%
_{n}(K)$. Then either $\Gamma $ contains a Zariski open solvable subgroup 
or $\Gamma $ contains a Zariski dense free subgroup of finite rank.
\end{thm}

\begin{rem}
Theorem \ref{Tits-alternative} seems quite close to the original theorem of Tits, stated above. 
And indeed, it is stated explicitly in \cite{Tits0} in the particular case when the Zariski closure of $\Gamma$ is assumed to be 
a semisimple Zariski connected algebraic group. However, the proof of Theorem \ref{Tits-alternative} relies 
on the methods developed in the present paper which allow one to deal with non Zariski connected groups. We will show below how Theorem \ref{Tits-alternative} can be easily deduced from Theorem \ref{main} below.
\end{rem}

The main purpose of our work is to prove the analog of Theorem \ref{Tits-alternative}, in case the ground 
field, and hence any linear group over it, carries a more interesting topology than the Zariski topology, namely for local fields.

Assume that $k$ is a local field, i.e. $\Bbb{R}$, $\Bbb{C}$, a finite extension
of $\Bbb{Q}_{p}$, or a field of formal power series in one variable over a
finite field. The full linear group GL$_{n}(k)$ and hence any subgroup of it, is endowed with the standard topology, that is the
topology induced from the local field $k.$ 
We then prove the following:

\begin{thm}[Topological Tits alternative]
\label{main} 
Let $k$ be a local field and $\Gamma $ a subgroup of GL$%
_{n}(K)$. Then either $\Gamma $ contains an open solvable subgroup 
or $\Gamma $ contains a dense free subgroup.
\end{thm}

Note that $\Gamma $ may contain both a dense free subgroup and an open
solvable subgroup: in this case $\Gamma $ has to be discrete and free. For
non discrete groups however, the two cases are mutually exclusive.

In general, the dense free subgroup from Theorem \ref{main} may have an
infinite (but countable) number of free generators. However, in many cases
we can find a dense free subgroup on finitely many free generators 
(see below Theorems \ref{precise0} and \ref{precise}). This is the case, for example, when $%
\Gamma $ itself is finitely generated. 
For another example consider the
group $\text{SL}_{n}({\mathbb{Q}})$, $n\geq2$. It is not finitely generated, yet, we
show that it contains a free subgroup of rank $2$ which is dense with
respect to the topology induced from $\text{SL}_{n}({\mathbb{R}})$.
Similarly, for any prime $p\in {\mathbb{N}}$, we show that $\text{SL}_{n}({%
\mathbb{Q}})$ contains a free subgroup of finite rank $r=r(n,p)\geq2$ which is
dense with respect to the topology induced from $\text{SL}_{n}({\mathbb{Q}}%
_{p})$.

When $char(k)=0$, the linearity assumption can be replaced by the weaker
assumption that $\Gamma $ is contained in some second-countable $k$-analytic
Lie group $G$. In particular, Theorem \ref{main} applies to subgroups of any
real Lie group with countably many connected components, and to subgroups of
any group containing a $p$-adic analytic pro-$p$ group as an open subgroup
of countable index.

Let us indicate how Theorem \ref{main} implies Theorem \ref{Tits-alternative}. Let $K$ be a field, $\gC\leq\GL_n(K)$ a finitely 
generated group, and let $R$ be the ring generated by the entries of $\gC$. By a variant of the Noether normalization theorem 
$R$ can be embedded in the valuation ring ${\mathcal{O}}$ of some local field $k$. Such
an embedding induces an embedding $i$ of $\Gamma $ in the linear pro-finite group $\text{GL}_{n}({\mathcal{O}})$. 
Note also that the topology induced on $\gC$ from the Zariski topology of $\GL_n(K)$ coincides with the one induced from the Zariski 
topology of $\GL_n(k)$ and this topology is weaker then the topology induced by the local field $k$. If $i(\gC )$ contains a relatively
open solvable subgroup then so is its closure, and by compactness, it follows that $\gC$ is almost solvable, and hence its Zariski 
connected component is solvable and Zariski open. If $i(\gC )$ does not contain an open solvable subgroup then, by Theorem \ref{main},
it contains a dense free subgroup which, as stated in a paragraph above, we may assume has finite rank. This free subgroup is indeed 
Zariski dense. 
   
The first step toward Theorem \ref{main} was carried out in our previous work \cite{BG}.
In \cite{BG} we made the assumption that $k=\BR$ and the closure of $\gC$ is connected.
This considerably simplifies the situation, mainly because it implies that $\gC$ is automatically 
Zariski connected. A main achievement of the present work is the understanding of dynamical properties 
of projective representations of non Zariski connected algebraic groups (see Section \ref{sec4}). 
Another  new aspect is the study of representations of finitely generated integral domains
into local fields (see Section \ref{Ti}) which allows us to avoid the rationality of the deformation space
of $\gC$ in $\GL_n(k)$, and hence to drop the assumption that $\gC$ is finitely generated. 

For the sake of simplicity, we restrict ourselves throughout this paper to a fixed local field. 
However, the proof of Theorem \ref{main} applies also in the following more general setup:

\begin{thm}
Let $k_1,k_2,\ldots ,k_r$ be local fields and let $\gC$ be a (finitely generated) subgroup of 
$\prod_{i=1}^r\GL_n(k_i)$. Assume that $\gC$
does not contain an open solvable subgroup, then $\gC$ contains a dense free subgroup (of finite rank).
\end{thm}

We also note that the argument of Section \ref{secinf}, where we build a dense free group on infinitely 
many generators, is applicable in a much grater generality. 
For example, we can prove the following adelic version:

\begin{prop}
Let $K$ be an algebraic number field and $\BG$ a simply connected semisimple algebraic group defined over $K$,
let $V_K$ be the set of all valuations of $K$. 
Then we have:
\begin{itemize}
\item {\it Strong approximation version:}
For any $v_0\in V_K$ such that $\BG$ in not $K_{v_0}$ anisotropic, $\BG (K)$ contains a free subgroup of 
infinite rank whose image under the diagonal embedding is dense in the restricted topological product 
corresponding to $V_K\setminus\{ v_0\}$.
\item {\it Weak approximation version:} 
$\BG (K)$ contains a free subgroup of infinite rank whose image under the diagonal embedding is dense 
in the direct product $\prod_{v\in V_K}\BG (K_v)$.
\end{itemize}
\end{prop}

Theorem \ref{main} has various applications. 
We shall now indicate some of them.

\subsection{Applications to the theory of pro-finite groups}

When $k$ is non-Archimedean, Theorem \ref{main} provides some new results
about pro-finite groups (see Section \ref{applications}). In particular, we
answer a conjecture of Dixon, Pyber, Seress and Shalev (cf. \cite{DPSS} and 
\cite{P}), by proving:

\begin{thm}\label{ppp}
Let $\Gamma $ be a finitely generated linear group over an arbitrary field.
Suppose that $\Gamma $ is not virtually solvable, then its pro-finite
completion $\hat{\Gamma}$ contains a dense free subgroup of finite rank.
\end{thm}

In \cite{DPSS}, using the classification of finite simple groups, the weaker
statement, that $\hat{\Gamma}$ contains a free subgroup whose closure is of
finite index, was established. Let us remark that the passage from a
subgroup whose closure is of finite index, to a dense subgroup is also a crucial
step in the proof of Theorem \ref{main}. It is exactly this problem that forces
us to deal with representations of non Zariski connected algebraic groups.
Additionally, our proof of \ref{ppp} does not rely on \cite{DPSS}, neither on the 
classification of finite simple groups.

We also note that $\Gamma $ itself
may not contain a pro-finitely dense free subgroup of finite rank. It was
shown in \cite{Sco} that surface groups have the property that any proper
finitely generated subgroup is contained in a proper subgroup of finite
index (see also \cite{SoVa}).

In Section \ref{applications} we also answer a conjecture of Shalev about
coset identities in pro-$p$ groups in the $p$-adic analytic case:

\begin{prop}
Let $G$ be an analytic pro-$p$ group. If $G$ satisfies a coset identity with respect to some open subgroup, 
then $G$ is virtually solvable, and in particular, satisfies an identity.
\end{prop}

\subsection{Applications in dynamics}

The question of existence of a free subgroup is closely related to
questions concerning amenability. It follows from Tits' alternative that for a finitely generated
linear group $\gC$, the following are equivalent:
\begin{itemize}
\item $\gC$ is amenable,
\item $\gC$ is almost solvable,
\item $\gC$ does not contain a non-abelian free subgroup.
\end{itemize}

The topology enters the game when considering actions of subgroups on
the full group. Let $k$ be a local field and $G\leq \GL_n(k)$ be a closed subgroup and $\gC\leq G$ a countable subgroup.
Let $P\leq G$ be any closed amenable subgroup, and consider the action of $\gC$ on $G/P$ by left 
multiplications.
Theorem \ref{main} implies:

\begin{thm}
The following are equivalent:\\
$(I)$ The action of $\gC$ on $G/P$ is amenable,\\
$(II)$ $\gC$ contains an open solvable subgroup,\\
$(III)$ $\gC$ does not contain a non-discrete free subgroup.
\end{thm} 

The equivalence between (I) and (II) for the Archimedean case (i.e. $k=\Bbb{R}$) was conjectured by Connes and Sullivan
and subsequently proved by Zimmer \cite{Zim} by means of super-rigidity methods.
The equivalence between (III) and (II) was asked by Carri\`{e}re and Ghys \cite{Ghys} who showed that (I) 
implies (III) (see also Section \ref{amenable}). 
For the case $G=\SL_2(\BR )$ they actually proved that (III) implies (II) and hence concluded the validity 
of the Connes-Sullivan conjecture for this case (before Zimmer). 
We remark that the proof of Carri\`{e}re and Ghys relies on the existence of an open subset of elliptic 
elements in $\SL_2(\BR)$ and hence does not apply to an arbitrary Lie group.

\begin{rem}
1. When $\gC$ is not both discrete and free, the conditions are also equivalent to: (III')
$\gC$ does not contain a dense free subgroup.\\
2. For $k$ Archimedean, (II) is equivalent to: (II') The connected component of the closure 
$\overline{\gC}^{\circ}$ is solvable.\\ 
3. The implication $(II)\rightarrow (III)$ is trivial and $(II)\rightarrow (I)$ follows easily from the 
definition of amenable actions.\\
\end{rem}

We also generalized Zimmer's theorem for
arbitrary locally compact groups as follows (see Section \ref{amenable}):

\begin{thm}
\label{111} Let $\Gamma $ be a countable subgroup of a locally compact
topological group $G$. Then the action of $\Gamma $ on $G$ (as well as on $%
G/P$ for $P\leq G$ closed amenable) by left multiplication is amenable, if
and only if $\Gamma $ contains a relatively open subgroup which is amenable
as an abstract group.
\end{thm}

As a consequence of Theorem \ref{111} we get the following generalization of Auslander's
theorem (see \cite{raghunathan} Theorem 8.24):

\begin{thm}
\label{112} Let $G$ be a locally compact topological group, let $P\leq G$ be
a closed normal amenable subgroup, and let $\pi :G\to G/P$ be the canonical
projection. Suppose that $H\leq G$ is a subgroup which contains a relatively
open amenable subgroup. Then $\pi (H)$ also contains a relatively open amenable subgroup.
\end{thm}

Theorem \ref{112} has many interesting conclusions. For example, it is well known that the original theorem 
of Auslander (Theorem \ref{112} for real Lie groups) directly implies Bieberbach's classical theorem that 
any compact Euclidean manifold is finitely covered by a torus (part of Hilbert's 18th problem).
As a consequence of the general Theorem \ref{112} we obtain some information on the structure of lattices in 
general locally compact groups. If $G=G_c\times G_d$ is a direct product of a connected semisimple Lie group 
and a locally compact totally disconnected group. Then, it is easy to see that, the projection of any lattice 
in $G$ to the connected factor lies between a lattice to its commensurator.
Such information is useful since it says (as follows from Margulis' commensurator criterion for arithmeticity)
that if this projection is not a lattice itself then it is a subgroup of the commensurator of some arithmetic 
lattice (which is, up to finite index, $G_c(\BQ )$). Theorem \ref{112} implies that similar statement holds for
general $G$ (see Proposition \ref{G_cxG_d}).

\subsection{The growth of leaves in Riemannian foliations}

Y. Carri\`{e}re's interest in the Connes-Sullivan conjecture stemmed from his study of 
the growth of leaves in Riemannian foliations. In \cite{Car1} Carri\`{e}re asked whether there is a dichotomy 
between polynomial and exponential growth. 
In order to study this problem, Carri\`{e}re defined the notion of \textit{local growth} for a subgroup of
a Lie group (see Definition \ref{local-growth}) and showed the equivalence of the growth type of a generic 
leaf and the local growth of the holonomy group of the foliation viewed as a subgroup of the corresponding structural Lie group associated to the Riemannian foliation (see 
\cite{Mol}).

Tits' alternative implies, with some additional argument for solvable non-nilpotent groups, 
the dichotomy between polynomial and exponential growth for finitely generated linear groups. 
Similarly, Theorem \ref{main}, with some additional argument based on its proof for solvable non-nilpotent 
groups, implies the analogous dichotomy for the local growth:

\begin{thm}\label{4444}
Let $\Gamma $ be a finitely generated dense subgroup of a
connected real Lie group $G$. If $G$ is nilpotent then $\Gamma $ has
polynomial local growth. If $G$ is not
nilpotent, then $\Gamma $ has exponential local growth.
\end{thm}

As a consequence of Theorem \ref{4444} we obtain:
 
\begin{thm}\label{55555}
Let $\mathcal{F}$ be a Riemannian foliation on a compact
manifold $M$. The leaves of $\mathcal{F}$ have polynomial growth if and only
if the structural Lie algebra of $\mathcal{F}$ is nilpotent. Otherwise, generic leaves have exponential 
growth.
\end{thm}

The first half of Theorem \ref{55555} was actually proved by Carri\`{e}re in 
\cite{Car1}. Using Zimmer's proof of the Connes-Sullivan conjecture, he
first reduced to the solvable case, then he proved the nilpotency of the
structural Lie algebra of $\mathcal{F}$ by a delicate direct argument (see
also \cite{Hae1}). He then asked whether the second half of this theorem
holds. Both parts of Theorem \ref{55555} follow from Theorem \ref{main} and
the methods developed in its proof. We remark that although the content of Theorem \ref{55555} is about
dense subgroups of connected Lie groups, its proof relies on methods developed in Section \ref{Ti} of the 
current paper. 

If we consider instead the growth of the holonomy cover of each leaf, then the dichotomy shown in Theorem \ref{55555} holds for every leaf. On the other hand, it is easy to give an example of a Riemannian foliation on a compact manifold in which the growth of a generic leaf is exponential while some of the leaves are compact (see below Section \ref{local-growth}).

\medskip

The strategy used in this article to prove Theorem \ref{main} consists in
perturbing the generators $\gamma _{i}$ of $\Gamma $ within $\Gamma $ and in
the topology of GL$_{n}(k),$ in order to obtain (under the assumption that $%
\Gamma $ has no solvable open subgroup) free generators of a free subgroup which is still dense in $\Gamma $. As it turns out, there exists an identity
neighborhood $U$ of some non virtually solvable subgroup $\Delta \leq \Gamma
,$ such that any selection of points $x_i$ in $U\gamma _{i}U$ generate a dense
subgroup in $\Gamma $. The argument used here to prove this claim depends on
whether $k$ is Archimedean, $p$-adic or of positive characteristic.

In order to find a free group, we use a variation of the ping-pong
method used by Tits, applied to a suitable action of $\gC$ on some projective space over some local field $f$. 
As in \cite{Tits0} the ping-pong players are the so-called proximal elements
(a proximal transformation is a transformation of $\Bbb{P}(f^{n})$ which
contracts almost all $\Bbb{P}(f^{n})$ into a small ball). 
However, the original method of Tits (via the use of high powers
of semisimple elements to produce ping-pong players) is not applicable to
our situation and a more careful study of the contraction properties of
projective transformations is necessary. 

An important difficulty lies in finding a representation $\rho $ of $\Gamma $ into some PGL$_{n}(f)$ for
some local field $f$ (which may or may not be isomorphic to $k$) such that the Zariski closure of 
$\rho (\Delta )$ acts strongly irreducibly (i.e. fixes no finite union of proper projective
subspaces) and such that $\rho (U)$ contains very proximal elements.
What makes this step much harder is the fact that $\gC$ may not be Zariski connected.
We handle this problem in Section \ref{sec4}.
We would like to note that we gain motivation and inspiration from the beautiful work of Margulis and Soifer
\cite{MS} where a similar difficulty arose. 

We then make use of the ideas developed in \cite{BG} and
inspired from \cite{AMS}, where it is shown how the dynamical properties of
a projective transformation can be read off on its Cartan decomposition.
This allows to produce a set of elements in $U$ which ``play ping-pong'' on the projective
space $\Bbb{P}(f^{n}),$ and hence generate a free group (see Theorem \ref
{free}). Theorem \ref{free} provides a very handy way to generate free
subgroups, as soon as some infinite subset of matrices with entries in a
given finitely generated ring (e.g. an infinite subset of a finitely
generated linear group) is given. 

The method used in \cite{Tits0} and in \cite{BG} 
to produce the representation $\rho$ is based on finding a representation of a finitely generated subgroup of
$\gC$ into $\GL_n(K)$ for some algebraic number field, and then to replace the 
number field by a suitable completion of it. However, in \cite{BG} and \cite{Tits0}, a lot of freedom was possible in the choice of $K$ and the representation 
into $\GL_n(K)$, and what played the main role there, was 
the appropriate choice of completion. This approach is no more applicable for the situation considered in 
this paper, and we are forced to choose both $K$ and the representation of $\gC$ in $\GL_n(K)$ in a careful 
way. For this purpose, we prove a result
(generalizing a lemma of Tits) asserting that in an arbitrary finitely
generated integral domain, any infinite set can be sent to an unbounded set
under an appropriate embedding of the ring into some local field (see
Section \ref{Ti}). This result is crucial in particular when dealing with
non finitely generated subgroups in Section \ref{secinf}. 
It is also used in the proof of the growth of leaves dichotomy, 
in Section \ref{Ti}. 
Our proof uses a striking simple fact,
originally due to P\'{o}lya in the case $k=\Bbb{C}$, about the inverse image
of the unit disc under polynomial transformations (see Lemma \ref{inv}).

Let us end this introduction by setting few notations that will be used throughout the
paper. The notation $H\leq G$ means that $H$ is a subgroup of the group $G$.
By $[G,G]\,$we denote the derived group of $G$, i.e. the group generated by
commutators. Given a group $\Gamma ,$ we denote by $d(\Gamma )\in {\mathbb{N}%
}$ the minimal size of a generating set of $\Gamma $. If $\Omega \subset G$
is a subset of $G$, then $\langle \Omega \rangle $ denotes the subgroup of $%
G $ generated by $\Omega $. If $\Gamma $ is a subgroup of an algebraic
group, we denote by $\overline{\Gamma }^{z}$ its Zariski closure. Note that
the Zariski topology on rational points does not depend on the field of
definition, that is if $V$ is an algebraic variety defined over a field $K$
and if $L$ is any extension of $K$, then the $K$-Zariski topology on $V(K)$
coincides with the trace of the $L$-Zariski topology on it. To avoid
confusion, we shall always add the prefix ``Zariski'' to any topological
notion regarding the Zariski topology (e.g. ``Zariski dense'', ``Zariski
open''). For the topology inherited from the local field $k$, however, we
shall plainly say ``dense'' or ``open'' without further notice (e.g. $\text{%
SL}_{n}({\mathbb{Z}})$ is open and Zariski dense in $\text{SL}_{n}({%
\mathbb{Z}}[1/p])$, where $k={\mathbb{Q}}_{p}$).

\setcounter{tocdepth}{1}
\tableofcontents


\section{A generalization of a lemma of Tits\label{Ti}}

In the original proof of the Tits alternative, Tits used an easy but crucial
lemma saying that given a finitely generated field $K$ and an element $%
\alpha \in K$ which is not a root of unity, there always is a local field $k$
and an embedding $f:K\to k$ such that $|f(\alpha )|>1$. A natural and useful
generalization of this statement is the following lemma:

\begin{lem}\label{Tit}\label{GTL}
Let $R$ be a finitely generated integral domain, and let $I\subset R$ be an infinite
subset. Then there exists a local field $k$ and an embedding $%
i:R\hookrightarrow k$ such that $i(I)$ is unbounded.
\end{lem}

As explained below, this lemma provides a straightforward way to build the
proximal elements needed in the construction of dense free subgroups.

Before giving the proof of Lemma \ref{GTL} let us point out a straightforward consequence:

\begin{cor}[Zimmer \cite{Zim3}, Theorems 6 and 7, and \cite{HaV} 6.26]
There is no faithful conformal action of an infinite Kazhdan group on the
Euclidean 2-sphere $S^{2}$.
\end{cor}

\begin{proof}
Suppose there is an infinite Kazhdan subgroup $\gC$ in $\SL_2(\BC )$, the group of conformal transformations of
$S^2$. Since $\gC$ has property (T), it is finitely generated, and hence, Lemma \ref{GTL} could be 
applied to yield a faithful 
representation of $\gC$ into $\SL_2(k)$ for some local field $k$, with unbounded image. However $\PSL_2(k)$ 
acts faithfully with compact isotropy groups by isometries on the hyperbolic space $\BH^3$ if $k$ is 
Archimedean, and on a tree if it is not. As $\gC$ has property-(T), it 
must fix a point (c.f. \cite{HaV} 6.4 and 6.23 or \cite{Zim3} Prop. 18) 
and hence lie in some compact group. A contradiction.
\end{proof}

When $R$ is integral over $\Bbb{Z}$, the lemma follows easily by considering
the diagonal embedding of $R$ into a product of finitely many completions of
its field of fractions. The main difficulty comes from the possible presence
of transcendental elements. Our proof of Lemma \ref{GTL} relies on the
following interesting fact. Let $k$ be a local field, and let $\mu =\mu _{k}$
denote the standard Haar measure on $k$, i.e. the Lebesgue measure if $k$ is
Archimedean, and the Haar measure giving measure $1$ to the ring of integers 
$\mathcal{O}_{k}$ of $k$ when $k$ is non-Archimedean. Given a polynomial $P$
in $k[X]$, let 
\begin{equation*}
A_{P}=\left\{ x\in k,|P(x)|\leq 1\right\} .
\end{equation*}

\begin{lem}
\label{inv}For any local field $k$, there is a constant $c=c(k)$ such that $%
\mu (A_{P})\leq c$ for any monic polynomial $P\in k[X]$.
\end{lem}

\proof%
Let $\overline{k}$ be an algebraic closure of $k$, and $P$ a monic
polynomial in $k[X]$. We can write $P(X)=\prod (X-x_{i})$ for some $x_{i}\in 
\overline{k}$. The absolute value of $k$ extends uniquely to an absolute
value in $\overline{k}$ (see \cite{Lang} XII, 4, Theorem 4.1 p. 482). Now if 
$x\in A_{P}$ then $|P(x)|\leq 1$, and hence 
\begin{equation*}
\sum \log \left| x-x_{i}\right| =\log |P(x)|\leq 0.
\end{equation*}
But $A_{P}$ is measurable and bounded, therefore, integrating with respect to 
$\mu ,$ we obtain 
\begin{equation*}
\sum \int_{A_{P}}\log \left| x-x_{i}\right| d\mu (x)=\int_{A_{P}}\sum \log
\left| x-x_{i}\right| d\mu (x)\leq 0.
\end{equation*}
The lemma will now follow from the following \textbf{claim}: \textit{for any
measurable set }$B\subset k$\textit{\ and any point }$z\in \overline{k}$, 
\begin{equation}
\int_{B}\log \left| x-z\right| d\mu (x)\geq \mu (B)-c,  \label{zououp}
\end{equation}
\textit{where }$c=c(k)>0$\textit{\ is some constant independent of }$z$%
\textit{\ and }$B$\textit{.}

Indeed, let $\widetilde{z}\in k$ be such that $|\widetilde{z}-z|=\min_{x\in
k}|x-z|$, then $|x-z|\geq |x-\widetilde{z}|$ for all $x\in k$, so 
\begin{equation*}
\int_{B}\log \left| x-z\right| d\mu (x)\geq \int_{B}\log \left| x-\widetilde{%
z}\right| d\mu (x)=\int_{B-\widetilde{z}}\log \left| x\right| d\mu (x).
\end{equation*}
Therefore, it suffices to show $(1)$ when $z=0$. But a direct computation
for each possible field $k$ shows that $-\int_{|x|\leq 1}\log |x|d\mu
(x)<\infty $. Therefore taking $c=\mu \left\{ x\in k,|x|\leq e\right\}
+|\int_{|x|\leq 1}\log |x|d\mu (x)|$ we obtain $(1)$. This concludes the
proof of the lemma.%
\endproof%

Lemma \ref{inv} was proved by P\'{o}lya in \cite{Po} for the case $k={%
\mathbb{C}}$ by means of potential theory. P\'{o}lya's proof gives the best
constant $c({\mathbb{C}})=\pi $. For $k={\mathbb{R}}$ one can show that the
best constant is $c({\mathbb{R}})=4$ and that it can be realized as the
limit of the sequence of lengths of the pre-image of $[-1,1]$ by the
Chebyshev polynomials (under an appropriate normalization of these
polynomials). In the real case, this result admits generalizations to
arbitrary smooth functions such as the Van der Corput lemma (see \cite{Car}
for a multi-dimensional analog). For $k$ non-Archimedean, the constant is
always $<2$ and it tends to $1$ as the residue field ${\mathcal{O}}_{k}/\pi {%
\mathcal{O}}_{k}$ gets larger.

Let us just explain how, with a little more consideration, one can improve
the constant $c$ in the above proof\footnote{%
Let us also remark that there is a natural generalization of Lemma \ref{inv}
to higher dimension which follows by an analogous argument: For any local
field $k$ and $n\in {\mathbb{N}}$, there is a constant $c(k,n)$, such that
for any finite set $\{ x_{1},\ldots ,x_{m}\}\in k^{n}$, we have $\mu \big(\{y\in
k^{n}:\prod_{1}^{m}\Vert y-x_{i}\Vert \leq 1\}\big)\leq c(k,n)$.}. We wish
to find the minimal $c>0$ such that for every compact subset $B$ of $k$
whose measure is $\mu (B)\geq c$ we have 
\begin{equation*}
\int_{B}\log |x|d\mu (x)\geq 0.
\end{equation*}
Suppose $k=\Bbb{C}$. Since $\log |x|$ is increasing with $|x|$, for any $B$ 
\begin{equation*}
\int_{B}\log |x|d\mu (x)\geq \int_{C}\log |x|d\mu (x)
\end{equation*}
where $C$ is a ball around $0$ ($C=\{x\in k:|x|\leq t\}$) with the same
area as $B$. Therefore $c=\pi t^{2}$ where $t$ is such that $2\pi
\int_{0}^{t}r\log (r)dr=0$. The unique positive root of this equation is $t=%
\sqrt{e}$. Thus we can take 
\begin{equation*}
c=\pi e.
\end{equation*}
For $k=\Bbb{R}$ the same argument gives a possible constant $c=2e$, while
for $k$ non-Archimedean it gives $c=1+\frac{1}{f(q-1)}$ where $q=p^{f}$ is
the size of the residue class field and $f$ is its dimension over its prime
field $\Bbb{F}_{p}$. 

As in the proof of Lemma \ref{inv}, there is a positive constant $c_{1}$
such that the integral of $\log |x|$ over a ball of measure $c_{1}$ centered
at $0$ is at least $1$. This implies:

\begin{cor}
\label{corll}For any monic polynomial $P\in k[X]$, the integral of $\log
|P(x)|$ over any set of measure grater then $c_{1}$ is at least the degree $%
d^{\circ }P$.
\end{cor}

We shall also need the following two propositions:

\begin{prop}
\label{Pr1} Let $k$ be a local field and $k_{0}$ its prime field. If $%
(P_{n})_{n}$ is a sequence of monic polynomials in $k[X]$ such that the
degrees $d^{\circ }P_{n}\to +\infty $ as $n\to \infty $, and $\xi
_{1},\ldots ,\xi _{m}$ are given numbers in $k$, then there exists a number $%
\xi \in k,$ transcendental over $k_{0}(\xi _{1},\ldots ,\xi _{m})$, such
that $(|P_{n}(\xi )|)_{n}$ is unbounded in $k$.
\end{prop}

\begin{proof}
Let $T$ be the set of numbers in $k$ which are transcendental over $%
k_{0}(\xi _{1},\ldots ,\xi _{m})$. Then $T$ has full measure. For every $r>0$
we consider the compact set 
\begin{equation*}
K_{r}=\left\{ x\in k:\forall n\text{ ~}|P_{n}(x)|\leq r\right\} .
\end{equation*}
We now proceed by contradiction. Suppose $T\subset \bigcup_{r>0}K_{r}$. Then
for some large $r$, we have $\mu (K_{r})\geq c_{1},$ where $c_{1}>0$ is the
constant from Corollary \ref{corll}. This implies

\begin{equation*}
d^{\circ }P_{n}\leq \int_{K_{r}}\log \big( |P_{n}(x)|\big) d\mu (x)\leq \mu
(K_{r})\log r,
\end{equation*}
contradicting the assumption of the proposition.
\end{proof}

\begin{prop}
\label{Pr2} If $(P_{n})_{n}$ is a sequence of distinct polynomials in $\Bbb{Z%
}[X_{1},\ldots ,X_{m}]$ such that $\sup_{n}d^{\circ }P_{n}<\infty $, then
there exist algebraically independent numbers $\xi _{1},\ldots ,\xi _{m}$ in 
$\Bbb{C}$ such that $(|P_{n}(\xi _{1},\ldots ,\xi _{m})|)_{n}$ is unbounded
in $\Bbb{C}$.
\end{prop}

\begin{proof}
Let $d=\max_{n}d^{\circ }P_{n}$ and let $T$ be the set of all $m$-tuples
of complex numbers algebraically independent over ${\mathbb{Z}}$. The $P_{n}$%
's lie in 
\begin{equation*}
\{P\in \Bbb{C}[X_{1},\ldots,X_{m}]:d^{\circ }P\leq d\}
\end{equation*}
which can be identified, since $T$ is dense and polynomials are continuous,
as a finite dimensional vector subspace $V$ of the $\Bbb{C}$-vector space of
all functions from $T$ to $\Bbb{C}$. Let $l=\dim _{\Bbb{C}}V$. Then, as it
is easy to see, there exist $(\overline{x}_{1},\ldots ,\overline{x}_{l})\in
T^{l},$ such that the evaluation map $P\mapsto \big( P(\overline{x}%
_{1}),\ldots ,P(\overline{x}_{l})\big)$ from $V$ to $\Bbb{C}^{l}$ is a
linear isomorphism from $V$ to $\BC^{\dim V}$. Since the $P_{n}$'s belong to a $\Bbb{Z}$-lattice in $V$%
, so does their image under the evaluation map. Since the $P_{n}$'s
are all distinct, $\{ P_{n}(\overline{x}_{i})\}$ is unbounded for an
appropriate $i\leq l.$ 
\end{proof}

\proof[Proof of Lemma \ref{GTL}]%
Let us first assume that the characteristic of the field of fractions of $R$
is $0$. By Noether's normalization theorem, $R\otimes _{\Bbb{Z}}\Bbb{Q}$ is
integral over $\Bbb{Q}[\xi _{1},\ldots ,\xi _{m}]$ for some algebraically
independent elements $\xi _{1},\ldots ,\xi _{m}$ in $R$. Since $R$ is
finitely generated, there exists an integer $l\in \Bbb{N}$ such that the
generators of $R$, hence all elements of $R$, are roots of monic polynomials
with coefficients in $S=\Bbb{Z}[\frac{1}{l},\xi _{1},\ldots ,\xi _{m}]$.
Hence $R_{0}:=R[\frac{1}{l}]$ is
integral over $S$. Let $F$ be the field of fractions of $R_{0}$ and $K$ that
of $S$. Then $F$ is a finite extension of $K$ and there are finitely many
embeddings $\sigma _{1},\ldots ,\sigma _{r}$ of $F$ into some (fixed) algebraic
closure $\overline{K}$ of $K$. Note that $S$ is integrally closed. Therefore
if $x\in R,$ the characteristic polynomial of $x$ over $F$ belongs to $S[X]$
and equals 
\begin{equation*}
\prod_{1\leq i\leq r}(X-\sigma _{i}(x))=X^{r}+\alpha _{r}(x)X^{r-1}+\ldots
+\alpha _{1}(x)
\end{equation*}
where each $\alpha _{i}(x)\in S$. Since $I$ is infinite, we can find $i_{0}$
such that $\{\alpha _{i_{0}}(x)\}_{x\in I}$ is infinite. This reduces the
problem to the case $R=S$, for if $S$ can be embedded in a local field $k$
such that $\{|\alpha _{i_{0}}(x)|\}_{x\in I}$ is unbounded, then for at
least one $i$, the $|\sigma _{i}(x)|$'s will be unbounded in some finite
extension of $k$ in which $F$ embeds.

So assume $I\subset S={\mathbb{Z}} [\frac{1}{l},\xi_1,\ldots,\xi_m]$ and
proceed by induction on the transcendence degree $m$.

The case $m=0$ is easy since $S=\Bbb{Z}[\frac{1}{l}]$ embeds discretely (by
the diagonal embedding) in the finite product $\BR\prod_{p\backslash l}\BQ_p$.

Now assume $m\geq 1$. Suppose first that the total degrees of the $x$'s in $%
I $ are unbounded. Then, for say $\xi _{m}$, $\sup_{x\in I}d_{\xi
_{m}}^{\circ }x=+\infty $. Let $a(x)$ be the dominant coefficient of $x$ in
its expansion as a polynomial in $\xi _{m}.$ Then $a(x)\in \Bbb{Z}[\frac{1}{l%
},\xi _{1},\ldots ,\xi _{m-1}]$ and is non zero. 

If $\{a(x)\}_{x\in I}$ is
infinite, then we can apply the induction hypothesis and find an embedding
of $\Bbb{Z}[\frac{1}{l},\xi _{1},\ldots ,\xi _{m-1}]$ into some local field $%
k$ for which $\{|a(x)|\}_{x\in I}$ is unbounded. Hence $I^{\prime }:=\{x\in
I,|a(x)|\geq 1\}$ is infinite. Now $\frac{x}{a(x)}$ is a monic polynomial in $%
k[\xi _{m}]$, so we can then apply Proposition \ref{Pr1} and extend the
embedding to ${\mathbb{Z}}[\frac{1}{l},\xi _{1},\ldots ,\xi _{m-1}][\xi
_{m}]=S$ in $k$, such that $\{\frac{x}{a(x)}\}_{x\in I^{\prime }}$ is unbounded in 
$k$. The image of $I$, under this embedding, is unbounded in $k$.

Suppose now that $\{a(x)\}_{x\in I}$ is finite. Then either $a(x)\in {%
\mathbb{Z}}[\frac{1}{l}]$ for all but finitely many $x$'s or not. In the
first case we can embed ${\mathbb{Z}}[\frac{1}{l},\xi _{1},\ldots ,\xi
_{m-1}]$ into either $\Bbb{R}$ or $\Bbb{Q}_{p}$ (for some prime $p$ dividing 
$l$) so that $|a(x)|\geq 1$ for infinitely many $x$'s, while in the second
case we can find $\xi _{1},\ldots ,\xi _{m-1}$ algebraically independent in $%
\Bbb{C}$, such that $|a(x)|\geq 1$ for infinitely many $x$ in $I$, Then, the
same argument as above, using Proposition \ref{Pr1} applies.

Now suppose that the total degrees of the $x$'s in $I$ are bounded. If for
some infinite subset of $I$, the powers of $\frac{1}{l}$ in the coefficients of $x$
(lying in $\Bbb{Z}[\frac{1}{l}]$) are bounded from above, then we can apply
Proposition \ref{Pr2} to conclude. If not, then for some prime factor $p$ of 
$l$, we can write $x=\frac{1}{p^{n(x)}}\tilde{x}$ where $\tilde{x}\in \Bbb{Z}%
_{p}[\xi _{1},\ldots ,\xi _{m}]$ with at least one coefficient of $p$-adic
absolute value $1$, and the $n(x)\in \Bbb{Z}$ are not bounded from above. By
compactness, we can pick a subsequence $(\tilde{x})_{x\in I^{\prime }}$
which converges in $\Bbb{Z}_{p}[\xi _{1},\ldots ,\xi _{m}]$, and we may
assume that $n(x)\to \infty $ on this subsequence. The limit will be a
non-zero polynomial $\tilde{x}_{0}$. Pick arbitrary algebraically
independent numbers $z_{1},\ldots ,z_{m}\in \Bbb{Q}_{p}$. The limit
polynomial $\tilde{x}_{0}$ evaluated at the point $(z_{1},\ldots ,z_{m})\in {%
\mathbb{Q}}_{p}^{m}$ is not $0$, and the sequence of polynomial $(\tilde{x}%
)_{x\in I^{\prime }}$ evaluated at $(z_{1},\ldots ,z_{m})$ tends to $\tilde{x%
}_{0}(z_{1},\ldots ,z_{m})\neq 0$. Hence $\big( x(z_{1},\ldots ,z_{m})\big)%
_{x\in I^{\prime }}$ tends to $\infty $ in ${\mathbb{Q}}_{p}$. Sending the $%
\xi _{i}$'s to the $z_{i}$'s we obtain the desired embedding. (Note that in
this case, after $p$ is selected, the specific values of the $z_{i}$'s are
not important.)

Finally, let us turn to the case when $char(k)=p>0$. The first part of the
argument remains valid : $R$ is integral over $S=\Bbb{F}_{q}[\xi _{1},\ldots
,\xi _{m}]$ where $\xi _{1},\ldots ,\xi _{m}$ are algebraically independent
over $\Bbb{F}_{q}$ and this enables to reduce to the case $R=S$. Then we
proceed by induction on the transcendence degree $m$. If $m=1$, then the
assignment $\xi _{1}\mapsto \frac{1}{t}$ gives the desired embedding of $S$
into $\Bbb{F}_{q}((t))$. Let $m\geq 2$ and note that the total degrees of
elements of $I$ are necessarily unbounded. From this point the proof works
verbatim as in the corresponding paragraphs above.%
\endproof%


\section{Contracting projective transformations\label{sec2}}

In this section and the next, unless otherwise stated, $k$ is assumed to be
a local field, with no assumption on the characteristic.

\subsection{Proximality and ping-pong\label{first}\label{ping}}

Let us first recall some basic facts about projective transformations on $%
\Bbb{P}(k^{n})$, where $k$ is a local field. For proofs and a detailed (and
self-contained) exposition, see \cite{BG}, Section 3. We let $\left\| \cdot
\right\| $ be the standard norm on $k^{n}$, i.e. the standard Euclidean norm
if $k$ is Archimedean and $\left\| x\right\| =\max_{1\leq i\leq n}|x_{i}|$
where $x=\sum x_{i}e_{i}$ when $k$ is non-Archimedean and $(e_{1},\ldots
,e_{n})$ is the canonical basis of $k^{n}$. This norm extends in the usual
way to $\Lambda ^{2}k^{n}$. Then we define the \textit{standard metric} on $%
\Bbb{P}(k^{n})$ by 
\begin{equation*}
d([v],[w])=\frac{\left\| v\wedge w\right\| }{\left\| v\right\| \left\|
w\right\| }.
\end{equation*}
With respect to this metric, every projective transformation is bi-Lipschitz
on $\Bbb{P}(k^{n})$. For $\epsilon \in (0,1)$, we call a projective
transformation $[g]\in $PGL$_{n}(k)$ \textbf{$\epsilon $-contracting} if
there exist a point $v_{g}\in {\mathbb{P}}^{n-1}(k),$ called an attracting
point of $[g],$ and a projective hyperplane $H_{g}$, called a repelling
hyperplane of $[g]$, such that $[g]$ maps the complement of the $\epsilon $%
-neighborhood of $H_{g}\subset \Bbb{P}(k^{n})$ (the repelling neighborhood
of $[g]$) into the $\epsilon $-ball around $v_{g}$ (the attracting
neighborhood of $[g]$). We say that $[g]$ is $\epsilon $\textbf{-very
contracting} if both $[g]$ and $[g^{-1}]$ are $\epsilon $-contracting. A
projective transformation $[g]\in $PGL$_{n}(k)$ is called \textbf{$%
(r,\epsilon )$-proximal} ($r>2\epsilon >0$) if it is $\epsilon $-contracting
with respect to some attracting point $v_{g}\in \Bbb{P}(k^{n})$ and some
repelling hyperplane $H_{g}$, such that $d(v_{g},H_{g})\geq r$. The
transformation $[g]$ is called \textbf{$(r,\epsilon )$-very proximal} if
both $[g]$ and $[g]^{-1}$ are $(r,\epsilon )$-proximal. Finally $[g]\,$is
simply called \textbf{proximal} (resp. \textbf{very proximal}) if it is $(r$%
\textbf{$,\epsilon )$}-proximal (resp. $(r$\textbf{$,\epsilon )$}-very
proximal) for some $r>2$\textbf{$\epsilon >0$. }

The attracting point $v_{g}$ and repelling hyperplane $H_{g}$ of an $%
\epsilon $-contracting transformation are not uniquely defined. Yet, if $[g]$
is proximal we have the following nice choice of $v_{g}$ and $H_{g}$.

\begin{lem}
\label{fix}Let $\epsilon \in (0,\frac{1}{4})$. There exist two constants $%
c_{1},c_{2}\geq 1$ (depending only on the local field $k$) such that if $[g]$
is an $(r,\epsilon )$-proximal transformation with $r\geq c_{1}\epsilon $
then it must fix a unique point $\overline{v}_{g}$ inside its attracting
neighborhood and a unique projective hyperplane $\overline{H}_{g}$ lying
inside its repelling neighborhood. Moreover, if $r\geq c_{1}\epsilon ^{2/3}$%
, then all positive powers $[g^{n}]$, $n\geq 1$, are $(r-2\epsilon
,(c_{2}\epsilon )^{\frac{n}{3}})$-proximal transformations with respect to
these same $\overline{v}_{g}$ and $\overline{H}_{g}$.
\end{lem}

Let us postpone the proof of this lemma till the next paragraph.

An $m$-tuple of projective transformations $a_{1},\ldots ,a_{m}$ is called a 
\textbf{ping-pong $m$-tuple} if all the $a_{i}$'s are $(r,\epsilon )$-very
proximal (for some $r>2\epsilon >0$) and the attracting points of $a_{i}$
and $a_{i}^{-1}$ are at least $r$-apart from the repelling hyperplanes of $%
a_{j}$ and $a_{j}^{-1}$, for any $i\neq j$. Ping-pong $m$-tuples give rise
to free groups by the following variant of the \textit{ping-pong lemma} (see 
\cite{Tits0} 1.1) :

\begin{lem}
If $a_{1},\ldots ,a_{m}\in $PGL$_{n}(k)$ form a ping-pong $m$-tuple, then $%
\langle a_{1},\ldots ,a_{m}\rangle $ is a free group of rank $m$.
\end{lem}

A \textit{finite} subset $F\subset $PGL$_{n}(k)$ is called $(m,r\mathbf{)}$%
\textbf{-separating} ($r>0$, $m\in \Bbb{N}$) if for every choice of $2m$
points $v_{1},\ldots ,v_{2m}$ in $\Bbb{P}(k^{n})$ and $2m$ projective
hyperplanes $H_{1},\ldots ,H_{2m}$ there exists $\gamma \in F$ such that 
\begin{equation*}
\min_{1\leq i,j\leq 2m}\{d(\gamma v_{i},H_{j}),d(\gamma
^{-1}v_{i},H_{j})\}>r.
\end{equation*}
A separating set and an $\epsilon $-contracting element for small $\epsilon $
are precisely the two ingredients needed to generate a ping-pong $m$-tuple.
This is summarized by the following proposition (see \cite{BG} Propositions
3.8 and 3.11).

\begin{prop}
\label{P-VP}\label{prop3}\label{prop0}Let $F$ be an $(m,r)$-separating set ($%
r<1$, $m\in \Bbb{N}$) in PGL$_{n}(k)$. Then there is $C\geq 1$ such that for
every $\epsilon $, $0<\epsilon <1/C,$ we have

$(i)$ If $[g]\in $PGL$_{n}(k)$ is an $\epsilon $-contracting transformation,
one can find an element $[f]\in F$, such that $[gfg^{-1}]$ is $C\epsilon $%
-very contracting.

$(ii)$ If $a_{1},\ldots ,a_{m}\in $PGL$_{n}(k)$, and $\gamma $ is an $%
\epsilon $-very contracting transformation, then there are $h_{1},\ldots
,h_{m}\in F$ and $g_{1},\ldots ,g_{m}\in F$ such that 
\begin{equation*}
(g_{1}\gamma a_{1}h_{1},g_{2}\gamma a_{2}h_{2},\ldots ,g_{m}\gamma
a_{m}h_{m})
\end{equation*}
forms a ping-pong $m$-tuple and hence are free generators of a free group.
\end{prop}


\subsection{The Cartan decomposition}\label{KAK}

Now let $\Bbb{H}$ be a Zariski connected reductive $k$-split algebraic $k$%
-group and $H=\Bbb{H}(k)$. Let $\Bbb{T}$ be a maximal $k$-split torus and $T=%
\Bbb{T}(k)$. Fix a system $\Phi $ of $k$-roots of $\Bbb{H}$ relative to $%
\Bbb{T}$ and a basis $\Delta $ of simple roots. Let $\Bbb{X}(\Bbb{T})$ be
the group of $k$-rational multiplicative characters of $\Bbb{T}$ and $%
V^{\prime }=\Bbb{X}(\Bbb{T})\otimes _{\Bbb{Z}}\Bbb{R}$ and $V$ the dual
vector space of $V^{\prime }$. We denote by $C^{+}$ the positive Weyl
chamber: 
\begin{equation*}
C^{+}=\left\{ v\in V:\forall \alpha \in \Delta ,~\alpha (v)>0\right\} .
\end{equation*}
The Weyl group will be denoted by $W$ and is identified with the quotient $%
N_{H}(T)/Z_{H}(T)$ of the normalizer by the centralizer of $T$ in $H$. Let $%
K $ be a maximal compact subgroup of $H$ such that $N_{K}(T)$ contains
representatives of every element of $W$. If $k$ is Archimedean, let $A$ be
the subset of $T\,$consisting of elements $t$ such that $|\alpha (t)|\geq 1$
for every simple root $\alpha \in \Delta $. And if $k$ is non-Archimedean,
let $A$ be the subset of $T$ consisting of elements such that $\alpha
(t)=\pi ^{-n_{\alpha }}$ for some $n_{\alpha }\in \Bbb{N}\cup \{0\}$ for any
simple root $\alpha \in \Delta $, where $\pi $ is a given uniformizer for $k$
(i.e. the valuation of $\pi $ is $1$). Then we have the following \textit{%
Cartan decomposition} (see Bruhat-Tits \cite{Tits}) 
\begin{equation}
H=KAK.  \label{Cartan}
\end{equation}
In this decomposition, the $A$ component is uniquely defined. We can
therefore associate to every element $g\in H$ a uniquely defined $a_{g}\in A$%
.

Then, in what follows,\textit{\ we define} $\chi (g)$ to be equal to $\chi (a_{g})$ for any
character $\chi \in \Bbb{X}(\Bbb{T})$ and element $g\in H$. Although this conflicts with the original meaning of $\chi(g)$ when $g$ belongs to the torus $\Bbb{T}(k)$, we will keep this notation throughout the paper. Thus we always
have $|\alpha (g)|\geq 1$ for any simple root $\alpha $ and $g\in H$.

Let us note that the above decomposition (\ref{Cartan}) is no longer true
when $\Bbb{H}$ is not assumed to be $k$-split (see Bruhat-Tits \cite{Tits}
or \cite{Qu} for the Cartan decomposition in the general case).

If $\Bbb{H=GL}_{n}$ and $\alpha $ is the simple root corresponding to the
difference of the first two eigenvalues $\lambda _{1}-\lambda _{2}$, then $%
a_{g}$ is a diagonal matrix $diag(a_{1}(g),\ldots,a_{n}(g))$ where $|\alpha
(g)|=|\frac{a_{1}(g)}{a_{2}(g)}|$. Then we have the following nice criterion
for $\epsilon $-contraction, which justifies the introduction of this notion
(see \cite{BG} Proposition 3.3).

\begin{lem}
\label{g contracts}\label{contract}Let $\epsilon <\frac{1}{4}$. If $|\frac{%
a_{1}(g)}{a_{2}(g)}|\geq 1/\epsilon ^{2}$, then $[g]\in $PGL$_{n}(k)$ is $%
\epsilon $-contracting on $\Bbb{P}(k^{n})$. Conversely, suppose $[g]$ is $%
\epsilon $-contracting on $\Bbb{P}(k^{n})$ and $k$ is non-Archimedean with
uniformizer $\pi $ (resp. Archimedean), then $|\frac{a_{1}(g)}{a_{2}(g)}%
|\geq \frac{\left| \pi \right| }{\epsilon ^{2}}$ (resp. $|\frac{a_{1}(g)}{%
a_{2}(g)}|\geq \frac{1}{4\epsilon ^{2}}$).
\end{lem}

The proof of Lemma \ref{fix}, as well as of Proposition \ref{prop3}, is
based on the latter characterization of $\epsilon $-contraction and on the
following crucial lemma (see \cite{BG} Lemmas 3.4 and 3.5) :

\begin{lem}
\label{open}\label{crit}Let $r,\epsilon \in (0,1]$. If $|\frac{a_{1}(g)}{%
a_{2}(g)}|\geq \frac{1}{\epsilon ^{2}}$, then $[g]$ is $\epsilon $-contracting with
respect to the repelling hyperplane 
\begin{equation*}
H_{g}=[span\{k^{\prime }{}^{-1}(e_{i})\}_{i=2}^{n}]
\end{equation*}
and the attracting point $v_{g}=[ke_{1}]$, where $g=ka_{g}k^{\prime }$ is a
Cartan decomposition of $g$. Moreover, $[g]$ is $\frac{\epsilon ^{2}}{r^{2}}$%
-Lipschitz outside the $r$-neighborhood of $H_{g}$. Conversely assume that
the restriction of $[g]$ to some open set $O\subset {\mathbb{P}}(k^{n})$ is $%
\epsilon $-Lipschitz, then $|\frac{a_{1}(g)}{a_{2}(g)}|\geq \frac{1}{2\epsilon}$.
\end{lem}


\subsection{The proof of Lemma \ref{fix}}

Given a projective transformation $[h]$ and $\delta >0$, we say that $(H,v)$
is a $\delta $\textit{-related pair} of a repelling hyperplane and
attracting point for $[h]$, if $[h]$ maps the complementary of the $\delta $%
-neighborhood of $H$ inside the $\delta $-ball around $v$.

The attracting point and repelling hyperplane of an $\delta$-contracting
transformation $[h]$ are not uniquely defined. However, note that if $\delta
<\frac{1}{4}$ then for any two $\delta$-related pairs of $[h]$ $%
(H_{h}^{i},v_{h}^{i}),~i=1,2$, we have $d(v_{h}^{1},v_{h}^{2})<2\delta$.
Indeed, since $\delta <\frac{1}{4}$, the union of the $\delta$
-neighborhoods of the $H_{h}^{i}$'s does not cover ${\mathbb{P}}(k^{n})$.
Let $p\in \Bbb{P}(k^{n})$ be a point lying outside this union, then $%
d([h]p,v_{h}^{i})<\delta$ for $i=1,2$.

Now consider two $\delta $-related pairs $(H_{h}^{i},v_{h}^{i}),~i=1,2$ of
some projective transformation $[h]$, satisfying $d(v_{h}^{1},H_{h}^{1})\geq
r$ and no further assumption on the pair $(H_{h}^{2},v_{h}^{2})$. Suppose
that $1\geq r>4\delta $. Then we claim that $\text{Hd}(H_{h}^{1},H_{h}^{2})%
\leq 2\delta $, where $\text{Hd}$ denotes the standard distance between
hyperplanes, i.e. the Hausdorff distance. (Note that $\text{Hd}%
(H^{1},H^{2})=max_{x\in H^{1}}\{\frac{|f_{2}(x)|}{\left\| x\right\| }\}$
where $f_{2} $ is the unique (up to sign) norm one functional whose kernel
is the hyperplane $H_{2}$ (for details see \cite{BG} section 3).) To see
this, notice that if $\text{Hd}(H_{h}^{1},H_{h}^{2})$ were greater than $%
2\delta $ then any projective hyperplane $H$ would contain a point outside
the $\delta $-neighborhood of either $H_{h}^{1}$ or $H_{h}^{2}$. Such a
point is mapped under $[h]$ to the $\delta $-ball around either $v_{h}^{1}$
or $v_{h}^{2}$, hence to the $3\delta $-ball around $v_{h}^{1}$. This in
particular applies to the hyperplane $[h^{-1}]H_{h}^{1}$. A contradiction to
the assumption $d(H_{h}^{1},v_{h}^{1})>4\delta $. We also conclude that when 
$r>8\delta ,$ then for any two $\delta $-related pairs $(H^{i},v^{i})$ $%
i=1,2 $ of $[h]$, we have $d(v^{i},H^{j})>\frac{r}{2}$ for all $i,j\in
\{1,2\}$.

Let us now fix an arbitrary $\epsilon $-related pair $(H,v)$ of the $%
(r,\epsilon )$-proximal transformation $[g]$ from the statement of Lemma \ref
{fix}. Let also $(H_{g},v_{g})$ be the hyperplane and point introduced in
Lemma \ref{open}. From Lemmas \ref{contract} and \ref{open}, we see that the
pair $(H_{g},v_{g})$ is a $C\epsilon $-related pair for $[g]$ for some
constant $C\geq 1$ depending only on $k$. Assume $d(v,H)\geq r>8C\epsilon $.
Then it follows from the above that the $\epsilon $-ball around $v$ is
mapped into itself under $[g]$, and that $d(v,H_{g})>\frac{r}{2}$. From
Lemma \ref{open}, we obtain that $[g]$ is $(\frac{4C\epsilon }{r})^{2}$%
-Lipschitz in this ball, and hence $[g^{n}]$ is $(\frac{4C\epsilon }{r}%
)^{2n} $-Lipschitz there. Hence $[g]$ has a unique fixed point $\overline{v}%
_{g}$ in this ball which is the desired attracting point for all the powers
of $[g] $. Note that $d(v,\overline{v}_{g})\leq \epsilon $.

Since $[g^{n}]$ is $(\frac{4C\epsilon }{r})^{2n}$-Lipschitz on some open
set, it follows from Lemma \ref{open} that $|\frac{a_{2}(g^{n})}{a_{1}(g^{n})%
}|\leq 2(\frac{4C\epsilon }{r})^{2n}$, and from Lemma \ref{contract} that $%
[g^{n}]$ is $2(\frac{4C\epsilon }{r})^{n}$-contracting. Moreover, it is now
easy to see that if $r>(4C)^{2}\epsilon $, then for every $2(\frac{%
4C\epsilon }{r})^{n}$-related pair $(H_{n},v_{n})$ for $[g^{n}]$ $n\geq 2$,
we have $d(\overline{v}_g,v_{n})\leq 4(\frac{4C\epsilon }{r})^{n}$. (To see
this apply $[g^{n}]$ to some point of the $\epsilon $-ball around $v$ which
lies outside the $2(\frac{4C\epsilon }{r})^{n}$-neighborhood of $H_{n}$).
Therefore $(H_{n},\overline{v}_g)$ is a $6(\frac{4C\epsilon }{r})^{n}$%
-related pair for $[g^{n}],$ $n\geq 2.$

We shall now show that the $\epsilon $-neighborhood of $H$ contains a unique 
$[g]$-invariant hyperplane which can be used as a common repelling
hyperplane for all the powers of $[g]$. The set $\mathcal{F}$ of all
projective points at distance at most $\epsilon $ from $H$ is mapped into
itself under $[g^{-1}]$. Similarly the set $\mathfrak{H}$ of all projective
hyperplanes which are contained in $\mathcal{F}$ is mapped into itself under 
$[g^{-1}]$. Both sets $\mathcal{F}$, and $\mathfrak{H}$ are compact with
respect to the corresponding Grassmann topologies. The intersection $%
\mathcal{F}_{\infty }=\cap [g^{-n}]\mathcal{F}$ is therefore non empty and
contains some hyperplane $\overline{H}_{g}$ which corresponds to any point
of the intersection $\cap [g^{-n}]\mathfrak{H}$. We claim that $\mathcal{F}%
_{\infty }=\overline{H}_{g}$. Indeed, the set $\mathcal{F}_{\infty }$ is
invariant under $[g^{-1}]$ and hence under $[g]$ and $[g^{n}]$. Since $%
(H_{n},\overline{v}_{g})$ is a $6(\frac{4C\epsilon }{r})^{n}$-related pair
for $[g^{n}],$ $n\geq 2$, and since $\overline{v}_g$ is ``far'' (at least $%
r-2\epsilon $ away) from the invariant set $\mathcal{F}_{\infty }$, it
follows that for large $n$, $\mathcal{F}_{\infty }$ must lie inside the $6(%
\frac{4C\epsilon }{r})^{n}$-neighborhood of $H_{n}$. Since $\mathcal{F}%
_{\infty }$ contains a hyperplane, and since it is arbitrarily close to a
hyperplane, it must coincide with a hyperplane. Hence $\mathcal{F}_{\infty }=%
\overline{H}_{g}$. It follows that $(\overline{H}_{g},\overline{v}_{g})$ is
a $12(\frac{4C\epsilon }{r})^{n}$-related pair for $[g^{n}]$ for any large
enough $n$. Note that then $d(\overline{v}_{g},\overline{H}_{g})>r-2\epsilon 
$, since $d(\overline{v}_{g},v)\leq \epsilon $ and $\text{Hd}(\overline{H}%
_{g},H)\leq \epsilon $. This proves existence and uniqueness of $(\overline{H%
}_{g},\overline{v}_{g})$ as soon as $r>c_{1}\epsilon $ where $c_{1}\geq
(4C)^{2}+8C$.

If we assume further that $r^{3}\geq 12(4C\epsilon )^{2}$, then $\mathcal{F}%
_{\infty }$ lies inside the $6(\frac{4C\epsilon }{r})^{n}$-neighborhood of $%
H_{n}$ as soon as $n\geq 2$. Then $(\overline{H}_{g},\overline{v}_{g})$ is a 
$12(\frac{4C\epsilon }{r})^{n}$-related pair for $[g^{n}]$, hence a $%
(c_{2}\epsilon )^{n/3}$-related pair for $[g^{n}]$ whenever $n\geq 1$, where 
$c_{2}\geq 1$ is a constant easily computable in terms of $C$. This finishes
the proof of the lemma. \endproof

In what follows, whenever we add the article \textit{the} to an attracting
point and repelling hyperplane of a proximal transformation $[g]$, we shall
mean these fixed point $\overline{v}_{g}$ and fixed hyperplane $\overline{H}%
_{g}$ obtained in Lemma \ref{fix}.


\subsection{The case of general semisimple group\label{gengen}}

Now let us assume that $\Bbb{H}$ is a Zariski connected semisimple $k$%
-algebraic group, and let $(\rho ,V_{\rho })$ be a finite dimensional $k$%
-rational representation of $\Bbb{H}$ with highest weight $\chi _{\rho }$.
Let $\Theta _{\rho }$ be the set of simple roots $\alpha $ such that $\chi
_{\rho }/\alpha $ is again a non-trivial weight of $\rho $ 
\begin{equation*}
\Theta _{\rho }=\left\{ \alpha \in \Delta :\chi _{\rho }/\alpha ~\text{is}~%
\text{a~weight~of}~\rho \right\} .
\end{equation*}
It turns out that $\Theta _{\rho }$ is precisely the set of simple roots $%
\alpha $ such that the associated fundamental weight $\pi _{\alpha }$
appears in the decomposition of $\chi _{\rho }$ as a sum of fundamental
weights. Suppose that the weight space $V_{\chi _{\rho }}$ corresponding to $%
\chi _{\rho }$ has dimension $1$, then we have the following lemma.

\begin{lem}
\label{second}There are positive constants $C_{1}\leq 1\leq C_{2}$, such
that for any $\epsilon \in (0,1)$ and any $g\in \Bbb{H}(k)$, if $|\alpha
(g)|>\frac{C_{2}}{\epsilon ^{2}}$ for all $\alpha \in \Theta _{\rho }$ then
the projective transformation $[\rho (g)]\in $PGL$(V_{\rho })$ is $\epsilon $%
-contracting, and conversely, if $[\rho (g)]$ is $\epsilon $-contracting,
then $|\alpha (g)|>\frac{C_{1}}{\epsilon ^{2}}$ for all $\alpha \in \Theta
_{\rho }$.
\end{lem}

\proof%
Let $V_{\rho }=\bigoplus V_{\chi }$ be the decomposition of $V_{\rho }$ into
a direct sum of weight spaces. Let us fix a basis $(e_{1},\ldots ,e_{n})$ of 
$V_{\rho }$ compatible with this decomposition and such that $V_{\chi _{\rho
}}=ke_{1}$. We then identify $V_{\rho }$ with $k^{n}$ via this choice of
basis. Let $g=k_{1}a_{g}k_{2}$ be a Cartan decomposition of $g$ in $H$. We
have $\rho (g)=\rho (k_{1})\rho (a_{g})\rho (k_{2})\in \rho (K)D\rho (K)$
where $D\subset SL_{n}(k)$ is the set of diagonal matrices. Since $\rho (K)$
is compact, there exists a positive constant $C$ such that if $[\rho (g)]$
is $\epsilon $-contracting then $[\rho (a_{g})]$ is $C\epsilon $%
-contracting, and conversely if $[\rho (a_{g})]$ is $\epsilon $-contracting
then $[\rho (g)]$ is $C\epsilon $-contracting. Therefore, it is equivalent
to prove the lemma for $\rho (a_{g})$ instead of $\rho (g)$. Now the
coefficient $|a_{1}(\rho (a_{g}))|$ in the Cartan decomposition on $%
SL_{n}(k) $ equals $\max_{\chi }|\chi (a_{g})|=|\chi _{\rho }(g)|$, and the
coefficient $|a_{2}(\rho (a_{g}))|$ is the second highest diagonal
coefficient and hence of the form $|\chi _{\rho }(a_{g})/\alpha (a_{g})|$
where $\alpha $ is some simple root. Now the conclusion follows from Lemma 
\ref{contract}. 
\endproof%


\section{Irreducible representations of non-Zariski connected algebraic
groups}\label{secsec}\label{sec4}

In the process of constructing dense free groups, we need to find some
suitable linear representation of the group $\Gamma $ we started with. In
general, the Zariski closure of $\Gamma $ may not be Zariski connected, and
yet we cannot pass to a subgroup of finite index in $\Gamma $ in Theorem \ref
{main}. Therefore we will need to consider representations of non Zariski
connected groups.

Let $\Bbb{H}^{\circ }$ be a connected semisimple $k$-split algebraic $k$%
-group. The group $Aut_{k}(\Bbb{H}^{\circ })$ of $k$-automorphisms of $\Bbb{H%
}^{\circ }$ acts naturally on the characters $\Bbb{X}(\Bbb{T})$ of a maximal
split torus $\Bbb{T}$. Indeed, for every $\sigma \in Aut_{k}(\Bbb{H}^{\circ
})$, the torus $\sigma (T)$ is conjugate to $T=\Bbb{T}(k)$ by some element $%
g\in H=\Bbb{H}(k)$ and we can define the character $\sigma (\chi )\,$by $%
\sigma (\chi )(t)=\chi (g^{-1}\sigma (t)g).$ This is not well defined, since
the choice of $g$ is not unique (it is up to multiplication by an element of
the normalizer $N_{H}(T)$). But if we require $\sigma (\chi )$ to lie in the
same Weyl chamber as $\chi $, then this determines $g$ up to multiplication
by an element from the centralizer $Z_{H}(T)$, hence it determines $\sigma
(\chi )$ uniquely. Note also that every $\sigma $ sends roots to roots and
simple roots to simple roots.

In fact, what we need are representations of algebraic groups whose
restriction to the connected component is irreducible. As explained below,
it turns out that an irreducible representation $\rho $ of a connected
semisimple algebraic group $\Bbb{H}^{\circ }$ extends to the full group $%
\Bbb{H}$ if and only if its highest weight is invariant under the action of $%
\Bbb{H}$ by conjugation.

We thus have to face the problem of finding elements in $\Bbb{H}^{\circ }(k)\cap \gC$
which are $\varepsilon $-contracting under such a representation $\rho $. By
Lemma \ref{second} this amounts to finding elements $h$ such that $\alpha
(h) $ is large for all simple roots $\alpha $ in the set $\Theta _{\rho }$
defined in Paragraph \ref{gengen}. As will be explained below, we can find
such a representation $\rho $ such that all simple roots belonging to $%
\Theta _{\rho }$ are images by some outer automorphisms $\sigma $'s of $\Bbb{%
H}^{\circ }$ (coming from conjugation by an element of $\Bbb{H}$) of a
single simple root $\alpha $. But $\sigma (\alpha )(h)$ and $\alpha (\sigma
(h))$ are comparable. The idea of the proof below is then to find elements $%
h $ in $\Bbb{H}^{\circ }(k)$ such that all relevant $\alpha (\sigma (h))$'s
are large. But, according to the converse statement in Lemma \ref{second},
this amounts to finding elements $h$ such that all relevant $\sigma (h)$'s
are $\varepsilon $-contracting under a representation $\rho _{\alpha }$ such
that $\Theta _{\rho _{\alpha }}=\{\alpha \}$. This is the content of the
forthcoming proposition.

Before stating the proposition, let us note that, $\Bbb{H}^{\circ}$ being $k$-split, to every simple root $%
\alpha \in \Delta $ corresponds an irreducible $k$-rational representation
of $\Bbb{H}^{\circ }(k)$ whose highest weight $\chi _{\rho _{\alpha }}$ is the fundamental weight $\pi _{\alpha }$ associated to $\alpha$ and has multiplicity one. In this case the set $\Theta _{\rho _{\alpha
}}$ defined in Paragraph \ref{gengen} is reduced to the singleton $\{\alpha
\}$.

\begin{prop}
\label{prop1}\label{prop2}Let $\alpha $ be a simple root. Let $I$ be a
subset of $\Bbb{H}^{\circ }(k)$ such that $\{|\alpha (g)|\}_{g\in I}$ is
unbounded in $\Bbb{R}$. Let $\Omega \subset \Bbb{H}^{\circ }(k)$ be a
Zariski dense subset. Let $\sigma _{1},\ldots ,\sigma _{m}$ be algebraic $k$%
-automorphisms of $\Bbb{H}^{\circ }$. Then for any arbitrary large $M>0$,
there exists an element $h\in \Bbb{H}^{\circ }(k)$ of the form $%
h=f_{1}\sigma _{1}^{-1}(g)\ldots f_{m}\sigma _{m}^{-1}(g)$ where $g\in I$
and the $f_{i}$'s belong to $\Omega $, such that $|\sigma _{i}(\alpha
)(h)|>M $ for all $1\leq i\leq m$.
\end{prop}

\proof%
Let $\epsilon \in (0,1)$ and $g\in I$ such that $|\alpha (g)|\geq \frac{1}{%
\epsilon ^{2}}$. Let $(\rho _{\alpha },V)$ be the irreducible representation
of $\Bbb{H}^{\circ }(k)$ corresponding to $\alpha $ as described above.
Consider the weight space decomposition $V_{\rho _{\alpha }}=\bigoplus
V_{\chi }$ and fix a basis $(e_{1},\ldots ,e_{n})$ of $V=V_{\rho _{\alpha }}$
compatible with this decomposition and such that $V_{\chi _{\rho _{\alpha
}}}=ke_{1}$. We then identify $V$ with $k^{n}$ via this choice of basis, and
in particular, endow $\Bbb{P}(V)$ with the standard metric defined in the
previous section. It follows from Lemma \ref{second} above that $[\rho
_{\alpha }(g)]$ is $\epsilon C$-contracting on $\Bbb{P}(V)$ for some
constant $C\geq 1$ depending only on $\rho _{\alpha }$. Now from Lemma \ref
{crit}, there exists for any $x\in \Bbb{H}^{\circ }(k)$ a point $u_{x}\in 
\Bbb{P}(V)$ such that $[\rho _{\alpha }(x)] $ is $2$-Lipschitz over some
open neighborhood of $u_{x}$. Similarly there exists a projective hyperplane 
$H_{x}$ such that $[\rho _{\alpha }(x)]$ is $\frac{1}{r^{2}}$-Lipschitz
outside the $r$-neighborhood of $H_{x}$. Moreover, combining Lemmas \ref
{contract} and \ref{crit} (and up to changing $C$ if necessary to a larger
constant depending this time only on $k$), we see that $[\rho _{\alpha }(g)]$
is $\frac{\epsilon ^{2}C^{2}}{r^{2}}$-Lipschitz outside the $r$-neighborhood
of the repelling hyperplane $H_{g}$ defined in Lemma \ref{crit}. We pick $%
u_{g}$ outside this $r$-neighborhood.

By modifying slightly the definition of a finite $(m,r)$-separating set (see
above Paragraph \ref{first}), we can say that a finite subset $F$ of $\Bbb{H}%
^{\circ }(k)$ is an $(m,r)$-separating set \textit{with respect to} $%
\rho_\ga $ and $\sigma _{1},\ldots ,\sigma _{m}$ if for every choice of $m$
points $v_{1},\ldots ,v_{m}$ in $\Bbb{P}(V)$ and $m$ projective hyperplanes $%
H_{1},\ldots ,H_{m}$ there exists $\gamma \in F$ such that 
\begin{equation*}
\min_{1\leq i,j,k\leq m,}d(\rho _{\alpha }(\sigma _{k}(\gamma
))v_{i},H_{j})>r>0.
\end{equation*}

\textbf{Claim : } The Zariski dense subset $\Omega $ contains a finite $%
(m,r) $-separating set with respect to $\rho _{\alpha }$ and $\sigma
_{1},\ldots ,\sigma _{m}$, for some positive number $r$.

\textit{Proof of claim : }For $\gamma \in \Omega ,$ we let $M_{\gamma }$ be
the set of all tuples $(v_{i},H_{i})_{1\leq i\leq m}$ such that there exists
some $i,j$ and $l$ for which $\rho _{\alpha }(\sigma _{l}(\gamma ))v_{i}\in
H_{j}$. Now $\bigcap_{\gamma \in \Omega }M_{\gamma }$ is empty, for
otherwise there would be points $v_{1},\ldots ,v_{m}$ in $P(V)$ and
projective hyperplanes $H_{1},\ldots ,H_{m}$ such that $\Omega $ is included
in the union of the closed algebraic $k$-subvarieties $\{x\in \Bbb{H}^{\circ
}(k) $, $\rho _{\alpha }(\sigma _{l}(x))v_{i}\in H_{j}\}$ where $i,j$ and $l$
range between $1$ and $m$. But, by irreducibility of $\rho _{\alpha }$ each
of these subvarieties is proper, and this would contradict the Zariski
density of $\Omega $ or the Zariski connectedness of $\Bbb{H}^{\circ }$.
Now, since each $M_{\gamma }$ is compact in the appropriate product of
Grassmannians, it follows that for some finite subset $F\subset \Omega $, $%
\bigcap_{\gamma \in F}M_{\gamma }=\emptyset $. Finally, since $\max_{\gamma
\in F}\min_{1\leq i,j,l\leq m}d(\rho _{\alpha }(\sigma _{l}(\gamma
)v_{i},H_{j})$ depends continuously on $(v_{i},H_{i})_{i=1}^{m}$ and never
vanishes, it must attain a positive minimum $r,$ by compactness of the set
of all tuples $(v_{i},H_{i})_{i=1}^{m}$ in 
\begin{equation*}
\big(\Bbb{P}(V)\times \Bbb{G}r_{\dim (V)-1}(V)\big)^{2m}.
\end{equation*}
Therefore $F$ is the desired $(m,r)$-separating set. \newline

Up to taking a bigger constant $C$, we can assume that $C$ is larger than
the bi-Lipschitz constant of every $\rho _{\alpha }(x)$ on $\Bbb{P}(k^{n})$
when $x$ ranges over the finite set $\{\sigma _{k}(f),f\in F,1\leq k\leq m\}$%
.

Now let us explain how to find the element $h=f_{m}\sigma _{m}^{-1}(g)\ldots
f_{1}\sigma _{1}^{-1}(g)$ we are looking for. We shall choose the $f_{j}$'s
recursively, starting from $j=1$, in such a way that all the elements $%
\sigma _{i}(h)$, $1\leq i\leq m$, will be contracting. Write 
\begin{eqnarray*}
\sigma _{i}(h) &=&\sigma _{i}(f_{m}\sigma _{m}^{-1}(g)\ldots f_{1}\sigma
_{1}^{-1}(g))= \\
\big(\sigma _{i}(f_{m})\sigma _{i}\sigma _{m}^{-1}(g)\cdot \ldots \cdot
\sigma _{i}(f_{i})\big) &g&\big( \sigma _{i}(f_{i-1})\sigma _{i}\sigma
_{i-1}(g)\cdot \ldots \cdot \sigma _{i}(f_{1})\sigma _{i}\sigma _{1}^{-1}(g)%
\big).
\end{eqnarray*}
In order to make $\sigma _{i}(h)$ contracting, we shall require that:

\begin{itemize}
\item  For $m\geq i\geq 2,$ $\sigma _{i}(f_{i-1})$ takes the image under $%
\sigma _{i}\sigma _{i-1}(g)\cdot \ldots \cdot \sigma _{i}(f_{1})\sigma
_{i}\sigma _{1}^{-1}(g)$ of some open set on which $\sigma _{i}\sigma
_{i-1}(g)\cdot \ldots \cdot \sigma _{i}(f_{1})\sigma _{i}\sigma _{1}^{-1}(g)$
is 2-Lipschitz, e.g. a small neighborhood of the point 
\begin{equation*}
u_{i}:=\big(\sigma _{i}\sigma _{i-1}(g)\cdot \ldots \cdot \sigma
_{i}(f_{1})\sigma _{i}\sigma _{1}^{-1}(g)\big) (u_{\sigma _{i}\sigma
_{i-1}(g)\cdot \ldots \cdot \sigma _{i}(f_{1})\sigma _{i}\sigma
_{1}^{-1}(g)})
\end{equation*}
at least $r$ apart from the hyperplane $H_{g}$, and

\item  For $m>j\geq i$, $\sigma _{i}(f_{j})$ takes the image of $u_{i}$
under $\big(\sigma _{i}\sigma _{j}^{-1}(g)\cdot \ldots \cdot \sigma
_{i}(f_{i})\big)g$ at least $r$ apart from the hyperplane $H_{\sigma
_{i}\sigma _{j+1}^{-1}(g)}$ of $\sigma _{i}\sigma _{j+1}^{-1}(g)$ (i.e. of
the next element on the left in the expression of $\sigma _{i}(h)$).
\end{itemize}

Assembling the conditions on each $f_{i}$ we see that there are $\leq m$
points that the $\sigma _{j}(f_{i})$'s$,~1\leq j\leq m$ should send $r$
apart from $\leq m$ projective hyperplanes.

This appropriate choice of $f_{1},\ldots ,f_{m}$ in $F$ forces each of $%
\sigma _{1}(h),\ldots ,\sigma _{m}(h)$ to be $\frac{2C^{m+2}\epsilon ^{2}}{%
r^{2m}}$-Lipschitz in some open subset of $\Bbb{P}(V)$. Lemma \ref{crit} now
implies that $\sigma _{1}(h),\ldots ,\sigma _{m}(h)$ are $C_{0}\epsilon $%
-contracting on $\Bbb{P}(V)$ for some constant $C_{0}$ depending only on $%
(\rho ,V)$.

Moreover $h\in Ka_{h}K$ and each of the $\sigma _{i}(K)$ is compact, we
conclude that $\sigma _{1}(a_{h}),\ldots ,\sigma _{m}(a_{h})$ are also $%
C_{1}\epsilon $-contracting on $\Bbb{P}(V)$ for some constant $C_{1}$. But
for every $\sigma _{i}$ there exists an element $b_{i}\in \Bbb{H}^{\circ
}(k) $ such that $\sigma _{i}(T)=b_{i}Tb_{i}^{-1}$ and $\sigma _{i}(\alpha
)(t)=\alpha (b_{i}^{-1}\sigma _{i}(t)b_{i})$ for every element $t$ in the
positive Weyl chamber of the maximal $k$-split torus $T=\Bbb{T}(k)$. Up to
taking a larger constant $C_{1}$ (depending on the $b_{i}$'s) we therefore
obtain that $b_{1}^{-1}\sigma _{1}(a_{h})b_{1},\ldots ,b_{m}^{-1}\sigma
_{m}(a_{h})b_{m}$ are also $C_{1}\epsilon $-contracting on $\Bbb{P}(V)$ via
the representation $\rho _{\alpha }$. Finally Lemma \ref{second} yields the
conclusion that $|\sigma _{i}(\alpha )(h)|=|\alpha (b_{i}^{-1}\sigma
_{i}(a_{h})b_{i})|\geq \frac{1}{C_{2}\epsilon ^{2}}$ for some other positive
constant $C_{2}$. Since $\epsilon $ can be chosen arbitrarily small, we are
done. 
\endproof%


Now let $\Bbb{H}$ be an arbitrary algebraic $k$-group, whose identity
connected component $\Bbb{H}^{\circ }$ is semisimple. Let us fix a system $%
\Sigma $ of $k$-roots for $\Bbb{H}^{\circ }$ and a simple root $\alpha $.
For every element $g$ in $\Bbb{H}(k)$ let $\sigma _{g}$ be the automorphism
of $\Bbb{H}^{\circ }(k)$ which is induced by $g$ under conjugation, and let $%
\mathcal{S}$ be the group of all such automorphisms. As was described above, 
$\mathcal{S}$ acts naturally on the set $\Delta $ of simple roots. Let $%
\mathcal{S}\cdot \alpha =\{\alpha _{1},\ldots ,\alpha _{p}\}$ be the orbit
of $\alpha $ under this action. Suppose $I\subset \Bbb{H}^{\circ }(k)$
satisfies the conclusion of the last proposition for $\mathcal{S}\cdot
\alpha $, that is for any $\epsilon >0$, there exists $g\in I$ such that $%
|\alpha _{i}(g)|>1/\epsilon ^{2}$ for all $i=1,\ldots ,p$. Then the
following proposition shows that under some suitable irreducible projective
representation of the full group $\Bbb{H}(k)$, for arbitrary small $\epsilon 
$, some elements of $I$ act as $\epsilon $-contracting transformations.

\begin{prop}
\label{prop4}Let $I\subset \Bbb{H}^{\circ }(k)$ be as above. Then there
exists a finite extension $K$ of $k,$ $[K:k]<\infty ,$ and a non-trivial
finite dimensional irreducible $K$-rational representation of $\Bbb{H}%
^{\circ }$ into a $K$-vector space $V$ which extends to an irreducible
projective representation $\rho :\Bbb{H}(K)\rightarrow $PGL$(V)$, satisfying
the following property : for every positive $\epsilon >0$ there exists $%
\gamma _{\epsilon }\in I$ such that $\rho (\gamma _{\epsilon })$ is an $%
\epsilon $-contracting projective transformation of $\Bbb{P}(V)$.
\end{prop}

\proof%
Up to taking a finite extension of $k$, we can assume that $\Bbb{H}^{\circ }$
is $k$-split. Let $(\rho ,V)$ be an irreducible $k$-rational representation
of $\Bbb{H}^{\circ }$ whose highest weight $\chi _{\rho }$ is a multiple of $%
\alpha _{1}+\ldots +\alpha _{p}$ and such that the highest weight space $%
V_{\chi _{\rho }}$ has dimension $1$ over $k$. Burnside's theorem implies
that, up to passing to a finite extension of $k,$ we can also assume that
the group algebra $k[\Bbb{H}^{\circ }(k)]$ is mapped under $\rho $ to the
full algebra of endomorphisms of $V$, i.e. $End_{k}(V)$. For a $k$%
-automorphism $\sigma $ of $\Bbb{H}^{\circ }$ let $\sigma (\rho )$ be the
representation of $\Bbb{H}^{\circ }$ given on $V$ by $\sigma (\rho )(g)=\rho
(\sigma (g))$. It is a $k$-rational irreducible representation of $\Bbb{H}%
^{\circ }$ whose highest weight is precisely $\sigma (\chi _{\rho })$. But $%
\chi _{\rho }=d(\alpha _{1}+\ldots +\alpha _{p})$ for some $d\in \Bbb{N}$,
and is invariant under the action of $\mathcal{S}.$ Hence for any $\sigma
\in \mathcal{S}$, $\sigma (\rho )$ is equivalent to $\rho $. So there must
exists a linear automorphism $J_{\sigma }\in $GL$(V)$ such that $\sigma
(\rho )(h)=J_{\sigma }\rho (h)J_{\sigma }^{-1}$ for all $h\in \Bbb{H}^{\circ
}(k)$. Now set $\widetilde{\rho }(g)=[\rho (g)]\in PGL(V)$ if $g\in \Bbb{H}%
^{\circ }(k)$ and $\widetilde{\rho }(g)=[J_{\sigma _{g}}]\in $PGL$(V)$
otherwise. Since the $\rho (g)$'s when $g$ ranges over $\Bbb{H}^{\circ }(k)$
generate the whole of $End_{k}(V)$, it follows from Schur's lemma that $%
\widetilde{\rho }$ is a well defined projective representation of the whole
of $\Bbb{H}(k)$. Now the set $\Theta _{\rho }$ of simple roots $\alpha $
such that $\chi _{\rho }/\alpha $ is a non-trivial weight of $\rho $ is
precisely $\{\alpha _{1},\ldots ,\alpha _{p}\}$. Hence if $\gamma _{\epsilon
}\in I\,$ satisfies $|\alpha _{i}(\gamma _{\epsilon })|>\frac{1}{\epsilon
^{2}}$ for all $i=1,\ldots ,p$, then we have by Lemma \ref{second} that $%
\widetilde{\rho }(\gamma _{\epsilon })$ is $C_2\epsilon$-contracting on $%
\Bbb{P}(V)$ for some constant $C_{2}$ independent of $\epsilon $. 
\endproof%

We can now state and prove the main result of this paragraph, and the only
one which will be used in the sequel. Let here $K$ be an arbitrary field which is finitely generated over its prime field and 
$\Bbb{H}$ an algebraic $K$-group such that its Zariski connected component $\Bbb{H}%
^{\circ }$ is semisimple and non-trivial. Fix some faithful $K$-rational representation 
$\Bbb{H}\hookrightarrow \Bbb{GL}_{d}$. Let $R$ be a finitely generated subring of $K$. We shall denote by $\Bbb{H}(R)$
(resp. $\Bbb{H}^{\circ }(R)$) the subset of points of $\Bbb{H}(K)$ (resp. $%
\Bbb{H}^{\circ }(K)$) which are mapped into $\Bbb{GL}_{d}(R)$ under the
latter embedding.

\begin{thm}\label{free} 
Let $\Omega _{0}\subset \Bbb{H%
}^{\circ }(R)$ be a Zariski dense subset of $\Bbb{H}^{\circ }$ with $\Omega
_{0}=\Omega _{0}^{-1}$. Suppose $\{g_{1},\ldots ,g_{m}\}$ is a finite subset
of $\Bbb{H}(K)$ exhausting all cosets of $\Bbb{H}^{\circ }$ in $\Bbb{H}$, and let 
$$
 \Omega =g_{1}\Omega _{0}g_{1}^{-1}\cup \ldots \cup g_{m}\Omega_{0}g_{m}^{-1}.
$$
Then we can find a number $r>0$, a local field $k$, an
embedding $K\hookrightarrow k,$ and a strongly irreducible projective representation 
$\rho :\Bbb{H}(k)\rightarrow PGL_{d}(k)$ defined over $k$ with the following
property. If $\epsilon \in (0,\frac{r}{2})$ and $a_{1},\ldots ,a_{n}\in \Bbb{%
H}(K)$ are $n$ arbitrary points ($n\in \Bbb{N}$), then there exist $n$
elements $x_{1},\ldots ,x_{n}$ with 
$$
 x_{i}\in \Omega ^{4m+2}a_{i}\Omega 
$$
such that the $\rho (x_{i})$'s form a ping-pong $n$-tuple of $(r,\epsilon )$%
-very proximal transformations on $\Bbb{P}(k^{d})$, and in particular are
generators of a free group $F_{n}$.
\end{thm}

\proof%
Up to enlarging the subring $R$ if necessary, we can assume that $K$ is the
field of fractions of $R$. We shall make use of Lemma \ref{Tit}. Since $%
\Omega _{0}$ is infinite, we can apply this lemma and obtain an embedding of 
$K$ into a local field $k$ such that $\Omega _{0}$ becomes an unbounded set
in $\Bbb{H}(k)$. Up to enlarging $k$ if necessary we can assume that $\Bbb{H}%
^{\circ }(k)$ is $k$-split. We fix a maximal $k$-split torus and a system of 
$k$-roots with a base $\Delta $ of simple roots. Then, in the corresponding
Cartan decomposition of $\Bbb{H}(k)$ the elements of $\Omega _{0}$ have
unbounded $A$ component (see Paragraph \ref{KAK}). Therefore, there exists a
simple root $\alpha $ such that the set $\{|\alpha (g)|\}_{g\in \Omega _{0}}$
is unbounded in $\Bbb{R}$. Let $\sigma _{g_{i}}$ be the automorphism of $%
\Bbb{H}^{\circ }(k)$ given by the conjugation by $g_{i}$. The orbit of $%
\alpha $ under the group generated by the $\sigma _{g_{i}}$'s is denoted by $%
\{\alpha _{1},\ldots ,\alpha _{p}\}$. Now it follows from Proposition \ref
{prop2} that for every $\epsilon >0$ there exists an element $h\in \Omega
^{2p}$ such that $|\alpha _{i}(h)|>1/\epsilon ^{2}$ for every $i=1,\ldots ,p$%
. We are now in a position to apply the last Proposition \ref{prop4} and
obtain (up to taking a finite extension of $k$ if necessary) an irreducible
projective representation $\rho :\Bbb{H}(k)\rightarrow $PGL$(V)$, such that
the restriction of $\rho $ to $\Bbb{H}^{\circ }(k)$ is also irreducible and
with the following property: for every positive $\epsilon >0$ there exists $%
h_{\epsilon }\in \Omega ^{2p}$ such that $\rho (h_{\epsilon })$ is an $%
\epsilon $-contracting projective transformation of $\Bbb{P}(V)$. Moreover,
since $\rho _{|\Bbb{H}^{\circ }}$ is also irreducible and $\Omega _{0}$ is
Zariski dense in $\Bbb{H}^{\circ }$ we can find an $(n,r)$-separating set
with respect to $\rho _{|\Bbb{H}^{\circ }}$ for some $r>0$ (for this
terminology, see definitions in Paragraph \ref{ping}). This follows from the
proof of the claim in Proposition \ref{prop1} above (see also Lemma 4.3. in 
\cite{BG}). By Proposition \ref{prop0} $(i)$ above, we obtain for every
small $\epsilon >0$ an $\epsilon $-very contracting element $\gamma
_{\epsilon }$ in $h_{\epsilon }\Omega _{0}h_{\epsilon }^{-1}\subset \Omega
^{4p+1}$. Similarly, statement $(ii)$ of the same Proposition gives elements 
$f_{1},\ldots ,f_{n}\in \Omega _{0}$ and $f_{1}^{\prime },\ldots
,f_{n}^{\prime }\in \Omega _{0}$ such that, for $\epsilon $ small enough, $%
(x_{1},\ldots ,x_{n})=(f_{1}^{\prime }\gamma _{\epsilon }a_{1}f_{1},\ldots
,f_{n}^{\prime }\gamma _{\epsilon }a_{n}f_{n})$ form under $\rho $ a
ping-pong $n$-tuple of proximal transformations on $\Bbb{P}(V)$. Then each $%
x_{i}$ lies in $\Omega ^{4p+2}a_{i}\Omega $ and together the $x_{i}$'s form
generators of a free group $F_{n}$ of rank $n$. 
\endproof%

\subsection{Further remarks}

For further use in later sections we shall state two more facts. Let $\Gamma
\subset \Bbb{G}(K)$ be a Zariski dense subgroup of some algebraic group $%
\Bbb{G}$. Suppose $\Gamma $ is not virtually solvable and let $\Delta \leq
\Gamma $ be a subgroup of finite index. Taking the quotient by the solvable
radical of $\Bbb{G}^{\circ }$, we obtain a homomorphism $\pi $ of $\Gamma $
into an algebraic group $\Bbb{H}$ whose connected component is semisimple.
Let $g_{1},\ldots ,g_{m}\in \Gamma $ be representatives of all different
cosets of $\Delta $ in $\Gamma $. Then $\Omega _{0}=\pi (\cap
_{i=1}^{m}g_{i}\Delta g_{i}^{-1})\cap \Bbb{H}^{0}$ is clearly Zariski dense
in $\Bbb{H}^{0}$ and satisfies the conditions of Theorem \ref{free}. Hence
taking $a_{i}=\pi (g_{i})$ in the theorem, we obtain:

\begin{cor}
\label{coset}\label{freecosets}Let $\Gamma $ be a linear group which is not
virtually solvable, and let $\Delta \subset \Gamma $ be a subgroup of finite
index. Then there is some choice of coset representatives for $\Gamma
/\Delta $ which generate a free group.
\end{cor}

The following lemma will be useful when dealing with the non-Archimedean
case.

\begin{lem}
\label{repgamma} Let $k$ be a non-Archimedean local field. Let $\Gamma \leq $%
GL$_{n}(k)$ be a linear group over $k$ which contains no open solvable
subgroup. Then there exists a homomorphism $\rho $ from $\Gamma $ into a $k$%
-algebraic group $\Bbb{H}$ such that the Zariski closure of the image of any
open subgroup of $\Gamma $ contains the connected component of identity ${%
\mathbb{H}}^{\circ }$. Moreover, we can take $\rho :\Gamma \rightarrow \Bbb{H%
}(k)$ to be continuous in the topology induced by $k,$ and we can find $\Bbb{%
H}$ such that $\Bbb{H}^{\circ }\,$is semisimple and $\dim ({\mathbb{H}}%
^{\circ })\leq \dim \overline{\Gamma }^{z}$.
\end{lem}

\proof%
Let $U_{i}$ be a decreasing sequence of open subgroups in $\text{GL}_{n}(k)$
forming a base of identity neighborhoods. Consider the decreasing sequence
of algebraic groups $\overline{\Gamma \cap U_{i}}^{z}$. This sequence must
stabilize after a finite step $s$. The limiting group $\Bbb{G\,}=\overline{%
\Gamma \cap U_{s}}^{z}$ must be Zariski connected. Indeed, the intersection
of $\Gamma \cap U_{s}$ with the Zariski connected component of identity of $%
\Bbb{G}$ is a relatively open subgroup and contains $\Gamma \cap U_{t}$ for
some large $t$. If $\Bbb{G}$ were not Zariski connected, then $\overline{%
\Gamma \cap U_{t}}^{z}$ would be a smaller algebraic group. Moreover, from
the assumption on $\Gamma $, we get that $\Bbb{G}$ is not solvable.

Note that the conjugation by an element of $\Gamma $ fixes $\Bbb{G}$, since $%
\gamma U_{i}\gamma ^{-1}\cap U_{i}\,$is again open if $\gamma \in \Gamma $
and hence contains some $U_{j}$. Since the solvable radical $Rad(\Bbb{G)}$
of $\Bbb{G}$ is a characteristic subgroup of $\Bbb{G}$, it is also fixed
under conjugation by elements of $\Gamma $. We thus obtain a homomorphism $%
\rho $ from $\Gamma $ to the $k$-points of the group of $k$-automorphisms $%
\Bbb{H=}Aut\Bbb{(S)}$ of the Zariski connected semisimple $k$-group $\Bbb{S=G%
}/Rad(\Bbb{G)}$. This homomorphism is clearly continuous. Since the image of 
$\Gamma \cap U$ is Zariski dense in $\Bbb{G}$ for all open $U\subset \text{GL%
}_{n}(k)$, $\Gamma \cap U$ is mapped under this homomorphism to a Zariski
dense subgroup of the group of inner automorphisms $Int(\Bbb{S})$ of $\Bbb{S}
$. But $Int(\Bbb{S})\,$is a semisimple algebraic $k $-group which is
precisely the Zariski connected component of identity of $\Bbb{H=}Aut\Bbb{(S)%
}$ (see for example \cite{Bo}, 14.9). Finally, it is clear from the
construction that $\dim {\mathbb{H}}\leq \dim \overline{\Gamma }^{z}$. 
\endproof%


\section{The proof of Theorem \ref{main} in the finitely generated case}

\label{proof}

In this section we prove our main result, Theorem \ref{main}, in the case
when $\Gamma $ is finitely generated. We obtain in fact a more precise
result which yields some control on the number of generators required for
the free group.

\begin{thm}
\label{precise0}Let $\Gamma \leq \text{GL}_{n}(k)$ be a finitely generated
linear group over a local field $k$. Suppose $\Gamma $ contains no solvable
open subgroup. Then, there is a constant $h(\Gamma )\in {\mathbb{N}}$ such
that for any integer $r\geq h(\Gamma )$, $\Gamma $ contains a dense free
subgroup of rank $r$. Moreover, if char$(k)=0$ we can take $h(\Gamma
)=d(\Gamma )$ (i.e. the minimal size of a generating set for $\Gamma $),
while if char$(k)>0$ we can take $h(\Gamma )=d(\Gamma )+n^{2}$.
\end{thm}

In the following paragraphs we split the proof to three cases (Archimedean,
non-Archimedean of characteristic zero, and positive characteristic) which
have to be dealt with independently.

\subsection{The Archimedean case\label{archi}}

Consider first the case $k={\mathbb{R}}$ or ${\mathbb{C}}$. Let $G$ be the
linear Lie group $G=\overline{\Gamma }$, and let $G^{\circ }$ be the
connected component of the identity in $G$. The condition ``$\Gamma $
contains no open solvable subgroup'' means simply ``$G^{\circ }$ is not
solvable''. Note also that $d(G/G^{\circ })\leq d(\Gamma )<\infty $.

Define inductively $G_{0}^{\circ }=G^{\circ }$ and $G_{n+1}^{\circ }=%
\overline{[G_{n}^{\circ },G_{n}^{\circ }]}$. This sequence stabilizes after
some finite step $t$ to a normal topologically perfect subgroup $%
H:=G_{t}^{\circ }$ (i.e. the commutator group $[H,H]\,$is dense in $H$). As
was shown in \cite{BG} Theorem 2.1, any topologically perfect group $H$
contains a finite set of elements $\{h_{1},\ldots ,h_{l}\},~l\leq \dim (H)$
and a relatively open identity neighborhood $V\subset H$ such that, for any
selection of points $x_{i}\in Vh_{i}V,$ the group $\langle x_{1},\ldots
,x_{l}\rangle $ is dense in $H$. Moreover, $H$ is clearly a characteristic
subgroup of $G^{\circ }$, hence it is normal in $G$. It is also clear from
the definition of $H$ that if $\Gamma $ is a dense subgroup of $G$ then $%
\Gamma \cap H$ is dense in $H$.

Let $r\geq d(\Gamma )$, and let $\{\gamma _{1},\ldots ,\gamma _{r}\}$ be a
generating set for $\Gamma $. Then one can find a smaller identity
neighborhood $U\subset V\subset H$ such that for any selection of points $%
y_{j}\in U\gamma _{j}U,~j=1,\ldots ,r,$ the group they generate $\langle
y_{1},\ldots ,y_{r}\rangle $ is dense in $G.$ Indeed, as $\Gamma \cap H$ is
dense in $H$, there are $l$ words $w_{i}$ in $r$ letters such that $%
w_{i}(\gamma _{1},\ldots ,\gamma _{r})\in Vh_{i}V$ for $i=1,\ldots ,l$.
Hence, for some smaller neighborhood $U\subset V\subset H$ and for any
selection of points $y_{j}\in U\gamma _{j}U,~j=1,\ldots ,r,$ we will have $%
w_{i}(y_{1},\ldots ,y_{r})\in Vh_{i}V$ for $i=1,\ldots ,l$. But then $%
\langle y_{1},\ldots ,y_{r}\rangle $ is dense in $G,$ since its intersection
with the normal subgroup $H$ is dense in $H$, and its projection to $G/H$
coincides with the projection of $\Gamma $ to $G/H$.

Let $R\leq G^{\circ }$ be the solvable radical of $G^{\circ }$. The group $%
G/R$ is a semisimple Lie group with connected component $G^{\circ }/R$ and $%
H $ clearly projects onto $G^{\circ }/R$. Composing the projection $G\to G/R$
with the adjoint representation of $G/R$ on its Lie algebra $\mathfrak{v}%
=Lie(G^{\circ }/R)$, we get a homomorphism $\pi :G\rightarrow $GL$(\mathfrak{%
v})$. The image $\pi (G)$ is open in the group of real points of some real algebraic
group ${\mathbb{H}}$ whose connected identity component $\Bbb{H}^{\circ }$
is semisimple. Moreover $\pi (\Gamma )$ is dense in $\pi (G)$. Let $m=|{%
\mathbb{H}}/{\mathbb{H}}^{\circ }|$ and let $g_{1},\ldots ,g_{m}\in \Gamma $
be elements which are sent under $\pi $ to representatives of all cosets of $%
{\mathbb{H}}^{\circ }\,$in $\Bbb{H}$. Let $U_{0}$ be an even smaller
symmetric identity neighborhood $U_{0}\subset U\subset H$ such that $%
U_{1}^{4m+2}\subset U\,$where $U_{1}=\cup _{j=1}^{m}g_{j}U_{0}g_{j}^{-1},$
and set $\Omega _{0}=\pi (U_{0}\cap \Gamma )$. Then the conditions of
Theorem \ref{free} are satisfied, since $\Omega _{0}$ is Zariski dense in $%
\Bbb{H}^{\circ }(\Bbb{R})$ (see \cite{BG} Lemma 5.2 applied to $H$). Thus we
can choose $\alpha _{i}\in \Gamma \cap U_{1}^{4m+2}\gamma _{i}U_{1}\,$which
generate a free group $\langle \alpha _{1},\ldots ,\alpha _{r}\rangle $. It
will also be dense by the discussion above.


\subsection{The $p$-adic case\label{padic}}

Suppose now that $k$ is a non-Archimedean local field of characteristic $0$,
i.e. it is a finite extension of the field of $p$-adic numbers ${\mathbb{Q}}%
_{p}$ for some prime $p\in {\mathbb{N}}$. Let ${\mathcal{O}}={\mathcal{O}}%
_{k}$ be the valuation ring of $k$ and $\frak{p}$ its maximal ideal. Let $%
\Gamma \leq \text{GL}_{n}(k)$ be a finitely generated linear group over $k$,
let $G$ be the closure of $\Gamma $ in GL$_{n}(k)$. Let $G({\mathcal{O}}%
)=G\cap $GL$_{n}({\mathcal{O}})$ (and $\Gamma ({\mathcal{O}})=\Gamma \cap G({%
\mathcal{O}})$) and denote by GL$_{n}^{1}({\mathcal{O}})$ the first
congruence subgroup, i.e. the kernel of the homomorphism GL$_{n}({\mathcal{O}%
})\rightarrow $GL$_{n}({\mathcal{O/}}\frak{p})$. The subgroup $G^{1}({%
\mathcal{O}})=G\cap \text{GL}_{n}^{1}({\mathcal{O}})$ is an open compact
subgroup of $G$ and is a $p$-adic analytic pro-$p$ group. The group GL$_{n}({%
\mathcal{O}})$ has finite rank (i.e. there is an upper bound on the minimal
number of topological generators for all closed subgroups of GL$_{n}({%
\mathcal{O}})$) as it follows for instance from Theorem 5.2 in \cite{Dix}.
Consequently, $G({\mathcal{O}})$ itself is finitely generated as a
pro-finite group and it contains the finitely generated pro-$p$ group $G^{1}(%
{\mathcal{O}})\,$as a subgroup of finite index. This implies that the
Frattini subgroup $F\leq G({\mathcal{O}})$ (the intersection of all maximal
open subgroups of $G({\mathcal{O}})$) is open (normal), hence of finite
index in $G({\mathcal{O}})\,$(Proposition 1.14 in \cite{Dix}). In this
situation, generating a dense group in $G$ is an open condition. More
precisely:

\begin{lem}
\label{3.2}Suppose $x_{1},\ldots ,x_{r}\in G$ generate a dense subgroup of $%
G $, then there is a neighborhood of identity $U\subset G$, such that for
any selection of points $y_{i}\in Ux_{i}U,~1\leq i\leq r$, the $y_{i}$'s
generate a dense subgroup of $G$.
\end{lem}

\proof
Note that a subgroup of the pro-finite group $G({\mathcal{O}})$ is dense if
and only if it intersects every coset of the Frattini subgroup $F$. Now
since $\langle x_{1},\ldots ,x_{r}\rangle $ is dense, there are $l=[G({%
\mathcal{O}}):F]$ words $\{w_{i}\}_{i=1}^{l}$ on $r$ letters, such that the $%
w_{i}(x_{1},\ldots ,x_{r})$'s are representatives of all cosets of $F$ in $G(%
{\mathcal{O}})$. But then, if $U$ is small enough, and $y_{i}\in Ux_{i}U$,
the elements $w_{i}(y_{1},\ldots ,y_{r})$ form again a full set of
representatives for the cosets of $F$ in $G({\mathcal{O}})$. This implies
that $\langle y_{1},\ldots ,y_{r}\rangle \cap G({\mathcal{O}})$ is dense in $%
G({\mathcal{O}})$. Now if we assume further that $U$ lies inside the open
subgroup $G({\mathcal{O}}),$ then we have $x_{i}\in G({\mathcal{O}})y_{i}G({%
\mathcal{O}})$, hence $x_{i}\in \overline{\langle y_{1},\ldots ,y_{r}\rangle 
}$ for $i=1,\ldots ,r$. This implies that $G=\overline{\langle y_{1},\ldots
,y_{r}\rangle }$. \endproof

The proof of the theorem now follows easily. Let $\{x_{1},\ldots ,x_{r}\}$
be a generating set for $\Gamma $. Choose $U$ as in the lemma, and take it
to be an open subgroup. Hence it satisfies $U^{l}=U$ for all $l\geq 1$. By
Lemma \ref{repgamma} we have a representation $\rho :\Gamma \to {\mathbb{H}}$
into some semisimple $k$-algebraic group ${\mathbb{H}}$ such that the image
of $U\cap \Gamma $ is Zariski dense in ${\mathbb{H}}^{0}$. Thus we can use
Theorem \ref{free} in order to find elements $\alpha _{i}\in \Gamma \cap
Ux_{i}U$ that generate a free group. It will be dense by Lemma \ref{3.2}.


\subsection{The positive characteristic case}

Finally, consider the case where $k$ is a local field of characteristic $p>0$%
, i.e. a field of formal power series ${\mathbb{F}}_{q}[[t]]$ over some
finite field extension ${\mathbb{F}}_{q}$ of ${\mathbb{F}}_{p}$. First, we
do not suppose that $\Gamma $ is finitely generated (in particular in the
lemma below). We use the same notations as those introduced at the beginning
of the last Paragraph \ref{padic} in the $p$-adic case. In particular $G$ is
the closure of $\Gamma $, $G(\mathcal{O})\,$is the intersection of $G$ with
GL$_{n}(\mathcal{O})$ where $\mathcal{O}$ is the valuation ring of $k$. In
positive characteristic, we have to deal with the additional difficulty
that, even when $\Gamma $ is finitely generated, $G(\mathcal{O})$ may not be
topologically finitely generated.

However, when $G(\mathcal{O})$ is topologically finitely generated, then the
argument used in the $p$-adic case (via Lemmas \ref{repgamma} and \ref{3.2})
applies here as well without changes. In particular, if $\overline{\Gamma }$
is compact, then we do not have to take more than $d(\Gamma )$ generators
for the dense free subgroup. This fact will be used in Section {\ref
{applications}}. We thus have:

\begin{prop}
\label{gC-compact} Let $k$ be a non-Archimedean local field and let ${%
\mathcal{O}}$ be its valuation ring. Let $\Gamma \leq \text{GL}_{n}({%
\mathcal{O}})$ be a finitely generated group which is not virtually
solvable, then $\Gamma $ contains a dense free group $F_{r}$ for any $r\geq
d(\Gamma )$.
\end{prop}

Moreover, it is shown in \cite{BaLa} that if $k$ is a local field of
positive characteristic, and $G={\mathbb{G}}(k)$ for some semisimple simply
connected $k$-algebraic group $\Bbb{G}$, then $G$ and $G({\mathcal{O}})$ are
finitely generated. Thus, the above proof applies also to this case and we
obtain:

\begin{prop}
\label{?} Let $k$ be a local field of positive characteristic, and let $G$
be the group of $k$ points of some semisimple simply connected $k$-algebraic
group. Let $\Gamma $ be a finitely generated dense subgroup of $G$, then $%
\Gamma $ contains a dense $F_{r}$ for any $r\geq d(\Gamma )$.
\end{prop}

\smallskip

Let us now turn to the general case, when $G(\mathcal{O})$ is not assumed
topologically finitely generated. As above, we denote by $\text{GL}_n^1(%
\mathcal{O})$ the first congruence subgroup $\text{Ker}\big(\text{GL}_{n}(%
\mathcal{O})\rightarrow \text{GL}_{n}(\mathcal{O}/\frak{p})\big)$. This
group is pro-$p$ and, as it is easy to see, the elements of torsion in GL$%
_{n}^{1}(\mathcal{O})\,$are precisely the unipotent matrices. In particular
the order of every torsion element is $\leq p^{n}$. Moreover, every open
subgroup of GL$_{n}^{1}(\mathcal{O})$ contains elements of infinite order.
Hence the torsion elements are not Zariski dense in GL$_{n}^{1}(\mathcal{O})$%
. More generally we have:

\begin{lem}
\label{tor}Let $k$ be a non-Archimedean local field of arbitrary
characteristic and $n$ a positive integer. There is an integer $m$ such that
the order of every torsion element in $\Bbb{GL}_{n}(k)$ divides $m$. In
particular, if $\Bbb{H}$ is a semisimple algebraic $k$-group, then the set
of torsion elements in $\Bbb{H}(k)$ is contained in a proper subvariety.
\end{lem}

\proof%
Let char$(k)=p\geq 0$. Suppose $x\in \Bbb{GL}_{n}(k)$ is an element of
torsion, then $x^{p^{n}}$ (resp. $x$ if $char(k)=0$) is semisimple. Since
the minimal polynomial of $x$ is of degree at most $n$, its eigenvalues,
which are roots of unity, lie in an extension of degree at most $n$ of $k$.
But, as $k$ is a local field, there are only finitely many such extensions.
Moreover, in a given non-Archimedean local field, there are only finitely
many roots of unity. Hence there is an integer $m$ such that $x^{m}=1$.

The last claim follows from the obvious fact that if $\Bbb{H}$ is semisimple
then there are elements if infinite order in $\Bbb{H}(k)$. 
\endproof%

Let $\rho :\Gamma \to {\mathbb{H}}$ be the representation given by Lemma \ref
{repgamma}. Then for any sufficiently small open subgroup $U\,$of GL$_{n}({%
\mathcal{O}})$ (for instance some small congruence subgroup), $\rho \big(%
\Gamma \cap U\big)$ is Zariski dense in ${\mathbb{H}}^{\circ }$. We then
have ($\Gamma $ is not assumed finitely generated):

\begin{lem}
\label{ff}There are $t:=\dim ({\mathbb{H}})$ elements $x_{1},\ldots
,x_{t}\in \Gamma ({\mathcal{O}})$ such that $\rho \big(\langle x_{1},\ldots
,x_{t}\rangle \big)$ is Zariski dense in ${\mathbb{H}}^{0}$.
\end{lem}

\proof
Let $U$ be an open subgroup of GL$_{n}({\mathcal{O}})\,$ so that $\rho \big(%
\Gamma \cap U\big)$ lies in $\Bbb{H}^{\circ }$ and is Zariski dense in it.
It follows from the above lemma that there is $x_{1}\in \Gamma \cap U$ such
that $\rho (x_{1})$ is of infinite order. Then the algebraic group $A=%
\overline{\langle x_{1}\rangle }^{z}$ is at least one dimensional.

Let the integer $i,$ $1\leq i\leq t,$ be maximal for the property that there
exist $i$ elements $x_{1},\ldots ,x_{i}\in \Gamma \cap U$ whose images in ${%
\mathbb{H}}$ generate a group whose Zariski closure is of dimension $\geq i$%
. We have to show that $i=t$. Suppose this is not the case. Fix such $%
x_{1},\ldots ,x_{i}$ and let $A$ be the Zariski connected component of
identity of $\overline{\langle \rho (x_{1}),\ldots ,\rho (x_{i})\rangle }%
^{z} $. Then for any $x\in \Gamma \cap U$, $\overline{\langle \rho
(x_{1}),\ldots ,\rho (x_{i}),\rho (x)\rangle }^{z}$ is $i$-dimensional. This
implies that $\rho (x)$ normalizes $A$. Since $\rho \big(\Gamma \cap U\big)$
is Zariski dense in ${\mathbb{H}}^{0}$, we see that $A$ is a normal subgroup
of ${\mathbb{H}}^{\circ }$. Dividing ${\mathbb{H}}^{\circ }$ by $A$ we
obtain a Zariski connected semisimple $k$-group of positive dimension, and a
map from $\Gamma \cap U$ with Zariski dense image into the $k$-points of
this semisimple group. But then, again by Lemma \ref{tor} above, there is an
element $x_{i+1}\in \Gamma \cap U$ whose image in ${\mathbb{H}}^{\circ }/A$
has infinite order --- a contradiction to the maximality of $i$. \endproof

Suppose now that $\Gamma $ is finitely generated and let $\Delta $ be the
closure in GL$_{n}(\mathcal{O})$ of the subgroup generated by $x_{1},\ldots
,x_{t}$ given by Lemma \ref{ff} above. It is a topologically finitely
generated pro-finite group containing the pro-$p$ subgroup of finite index $%
\Delta \cap G^{1}(\mathcal{O})$. Hence its Frattini subgroup $F$ is open and
of finite index (\cite{Dix} Proposition 1.14). In particular, $\rho (F\cap
\langle x_{1},\ldots ,x_{t}\rangle )$ is Zariski dense in ${\mathbb{H}}%
^{\circ }$, and we can use Theorem \ref{free} with $\Omega _{0}=\rho (F\cap
\langle x_{1},\ldots ,x_{t}\rangle )$. Note that $F$, being a group,
satisfies $F^{m}=F$ for $m\in {\mathbb{N}}$. Also $F$ is normal in $\Delta $%
. Let $\gamma _{1},\ldots ,\gamma _{r}$ be generators for $\Gamma $. By
Theorem \ref{free}, we can choose $\alpha _{i}\in F\gamma _{i}F,~i=1,\ldots
,r$, and $\alpha _{i+r}\in x_{i}F,~i=1,\ldots ,t$, so that $D=\langle \alpha
_{1},\ldots ,\alpha _{r+t}\rangle \,$is isomorphic to the free group $%
F_{r+t} $ on $r+t$ generators. Clearly $D\cap \Delta $ is dense in $\Delta $
and $D\cap F$ is dense in $F$. This implies that each $\gamma _{i}$ lies in $%
F\alpha _{i}F\subset \overline{D}$. As the $\gamma _{i}$'s generate $\Gamma
, $ we see that $D$ is dense in $\overline{\Gamma }\,$and this finishes the
proof.


\subsection{A stronger statement}

The argument above combined with the argument of \cite{BG} Section 2
provides the following generalization:

\begin{thm}
\label{precise}\label{??} Let $k$ be a local field and let $G\leq \text{GL}%
_{n}(k)$ be a closed linear group containing no open solvable subgroup.
Assume also that $G({\mathcal{O}})$ is topologically finitely generated in
case $k$ is non-Archimedean of positive characteristic. Then there is an
integer $h(G)$ which satisfies

\begin{itemize}
\item  $h(G)\leq 2\dim (G)-1+d(G/G^{\circ })$ if $k$ is Archimedean, and

\item  $h(G)$ is the minimal cardinality of a set generating a dense
subgroup of $G$ if $k$ is non-Archimedean,
\end{itemize}

such that any finitely generated dense subgroup $\Gamma \leq G$ contains a
dense $F_{r}$, for any $r\geq \min \{d(\Gamma ),h(G)\}$. Furthermore, if $k$
is non-Archimedean or if $G^{\circ }$ is topologically perfect (i.e. $%
\overline{[G^{\circ },G^{\circ }]}=G^{\circ }$) and $d(G/G^{\circ })<\infty $%
, then we can drop the assumption that $\Gamma $ is finitely generated. In
these cases, any dense subgroup $\Gamma $ in $G$ contains a dense $F_{r}$
for any $r\geq h(G)$.
\end{thm}

\begin{rem}
The interested reader is referred to \cite{BG} for a sharper estimation of $%
h(G)$ in the Archimedean case. For instance, if $G$ is a connected and
semisimple real Lie group, then $h(G)=2$.
\end{rem}

\medskip

Let us also remark that in the characteristic zero case, we can drop the
linearity assumption, and assume only that $\Gamma $ is a subgroup of some
second countable $k$-analytic Lie group. To see this, simply note that the
procedure of generating a dense subgroup does not rely upon the linearity of 
$G=\overline{\Gamma }$, and for generating a free subgroup, we can look at
the image of $G$ under the adjoint representation which is a linear group.
The main difference in the positive characteristic case is that we do not
know in that case whether or not the image $\text{Ad}(G)$ is solvable. For
this reason we make the additional linearity assumption in positive
characteristic.


\section{Dense free subgroups with infinitely many generators\label{secinf}}

In this section we shall prove the following:

\begin{thm}
\label{infinite} 
Let $k$ be a local field and $\gC\leq\GL_n(k)$ a linear group over $k$.
Assume that $\Gamma $ contains no open solvable subgroup,
then $\gC$ contains a countable dense free subgroup of infinite rank.
\end{thm}

As above, we denote by $G$ the closure of $\Gamma $. Since $\GL_n(k)$ and hence $G$ is second countable,
we can assume that $%
\Gamma$ is countable. If $k$ is non-Archimedean we let ${\mathcal{O}}$
denote the valuation ring of $k$, $G({\mathcal{O}}):=G\cap \text{GL}_{n}(%
{\mathcal{O}})$ the corresponding open pro-finite group and $\Gamma ({%
\mathcal{O}})=\Gamma \cap \text{GL}_{n}({\mathcal{O}})$. Set $G_{j}:=G\cap $%
GL$_{n}^{j}({\mathcal{O}})$ where GL$_{n}^{j}({\mathcal{O}})=\text{Ker}\big(%
\text{GL}_{n}(\mathcal{O})\rightarrow \text{GL}_{n}( \mathcal{O}/\frak{p}%
^{j})\big)$ is the $j$'th congruence subgroup, and write $\Gamma
_{j}:=\Gamma \cap G_{j}$. In order to treat both the Archimedean and the
non-Archimedean cases at the same time, we will say \textit{by convention}
that in the Archimedean case, $\Gamma _{j}\,$denotes always the same group $%
H\cap \Gamma $ where $H$ is the limit of the sequence of closed commutators
introduced in Paragraph \ref{archi}.

We now fix once and for all a sequence $(x_{j})$ of elements of $\Gamma $
which is dense in $G$. In the non-Archimedean case, we can also require that
the $G_{k_{j}}x_{j}G_{k_{j}}$'s form a base for the topology of $G$ for some
choice of a sequence of integers $(k_{j})$. We are going to perturb the $%
x_{i}$'s by choosing elements $y_{i}$'s inside $\Gamma_{k_{j}}x_{j}%
\Gamma_{k_{j}}$ which play ping-pong all together on some projective
space, and hence generate a dense free group.

From Lemma \ref{repgamma} in the non-Archimedean case, and from the
discussion in Paragraph \ref{archi} in the Archimedean case, we have a
homomorphism $\pi :\Gamma\rightarrow\Bbb{H}(k),$ where $\Bbb{H}$ is an
algebraic $k$-group with $\Bbb{H}^{\circ }$ semisimple, such that the
Zariski closure of $\pi (\Gamma _{j})$ contains $\Bbb{H}^{\circ }\,$for all $%
j\geq 1$. It now follows from Lemma \ref{ff} when $char(k)>0$ and from the
discussion in Paragraphs \ref{archi} and \ref{padic} in the other cases
(i.e. from the fact that $H$ and $G_{j}$ are topologically finitely
generated) that $\Gamma _{1}$ contains a finitely generated subgroup $\Delta
_{1}$ such that $\pi (\Delta _{1})$ is also Zariski dense in $\Bbb{H}^{\circ
}$ (we also take $\Delta_1$ to be dense in $H$ when $k$ is Archimedean).
From Theorem \ref{free} we can find a local field $k^{\prime }$ and an
irreducible projective representation $\rho$ of $\Bbb{H}$ on $\Bbb{P}(V_{k^{\prime }})$ 
defined over $k^{\prime }$ such that, under this representation, some
elements of $\Delta _{1}$ play ping-pong in the projective space $\Bbb{P}%
(V_{k^{\prime }}).$ In particular, for some $r>0$ and for every positive $%
\epsilon <\frac{r}{2}$, there is an element in $\Delta _{1}$ acting on $\Bbb{%
P}(V_{k^{\prime }})\,$by an $(r,\epsilon )$-very proximal transformation
(c.f. Paragraph \ref{ping}). Furthermore, there is a field extension $K$ of $%
k^{\prime }$ such that under this representation the full group $\Gamma $ is
map into PGL$(V_{K})$ where $V_{K}=V_{k^{\prime }}\otimes K$. This field
extension may not be finitely generated. Nevertheless, the
absolute value on $k^{\prime }\,$extends to an absolute value on $K$
(see \cite{Lang} XII, 4, Theorem 4.1 p. 482) and the projective space $\Bbb{P%
}(V_{K})\,$is still a metric space (although not compact in general) for the
metric introduced in Paragraph \ref{ping}. Moreover, if $[g]\in $PGL$%
(V_{k^{\prime }})\,$ is $\epsilon $-contracting on $\Bbb{P}(V_{k^{\prime }})$%
, it is $c\epsilon $-contracting on $\Bbb{P}(V_{K})$ for some constant $%
c=c(k^{\prime},K)\geq 1$. Similarly, if $[g]\in\text{PGL}_n(V_{k^{\prime}})$
is $(r,\epsilon )$-proximal transformation on ${\mathbb{P}} (V_{k^{\prime}})$
then it is $(\frac{r}{c},c\epsilon )$-proximal on ${\mathbb{P}} (V_{K})$.
Let $\rho :\Gamma \rightarrow $PGL$(V_{K})$ be this representation. (The
reason why we may not reduce to the case where $K$ is local is that there may not be a
finitely generated dense subgroup in $\Gamma $.)

In the Archimedean case, the discussion in Paragraph \ref{archi} shows that
we can find inside $\Delta_{1}$ elements $z_{1},\ldots ,z_{l}$, generating a
dense subgroup of $\Gamma _{1}=H\cap\Gamma$, and such that, under the above
representation, they act as a ping-pong $l$-tuple of projective
transformations. We can find another element $g\in \Delta_{1}\,$such that $%
(z_{1},\ldots,z_{l},g)$ acts as a ping-pong $(l+1)$-tuple and $g$ acts as an 
$(r,\epsilon )$-very proximal transformation on $\Bbb{P}(V_{k^{\prime }})$
where the pair $(r,\epsilon )$ satisfies the conditions of Lemma \ref{fix}
with respect to $k^{\prime}$. In the non-Archimedean case, let simply $g$ be
some element of $\Delta_{1}$ acting as an $(r,\epsilon )$-very proximal
transformation on the projective space $\Bbb{P}(V_{k^{\prime }})$ with $%
(r,\epsilon )$ as in Lemma \ref{fix}. As follows from Lemma \ref{fix}, $g$
(resp. $g^{-1}$) fixes an attracting point $\overline{v}_{g}$ (resp. $%
\overline{v}_{g^{-1}}$) and a repelling hyperplane $\overline{H}_{g}$ (resp. 
$\overline{H}_{g^{-1}}$) and the positive (resp. negative) powers $g^{n}$
behave as $(\frac{r}{C},(C\epsilon )^{\frac{n}{3}})$-very proximal
transformations with respect to these same attracting points and repelling
hyperplanes. Note that in the non-Archimedean case, if $n_{j}$ is the index
of the $j$'th congruence subgroup $G_{j}$ in $G({\mathcal{O}})$, then $%
g^{n_{j}}\in G_{j}$ and in particular $g^{n_{j}}\to 1$ as $j$ tends to
infinity.

We are now going to construct an infinite sequence $(g_{j})$ of elements in $%
\Delta _{1}$ acting on $\Bbb{P}(V_{k^{\prime }})$ by very proximal
transformations and such that they play ping-pong all together on $\Bbb{P}%
(V_{k^{\prime }})$ (and also together with $z_{1},\ldots,z_{l}$ in the
Archimedean case). Since $\pi (\Delta _{1})$ is Zariski dense in $\Bbb{H}%
^{\circ }$ and the representation $\rho$ of $\Bbb{H}^{\circ }\,$ is irreducible, we
may pick an element $\gamma\in \Delta _{1}$ such that 
\begin{equation*}
\{\rho (\gamma)\overline{v}_{g},~\rho (\gamma)\overline{v}_{g^{-1}},~\rho
(\gamma^{-1})\overline{v}_g, ~\rho (\gamma^{-1})\overline{v}_{g^{-1}}\}\cap %
\big(\overline{H}_{g}\cup\overline{H}_{g^{-1}}\cup\{\overline{v}_{g},%
\overline{v}_{g^{-1}}\}\big)=\emptyset.
\end{equation*}

Now consider the element $\delta_{m_1}=g^{m_1}\gamma g^{m_1}$. When ${m_1}$
is large enough, $\delta_{m_1}$ acts on ${\mathbb{P}} (V_{k^{\prime}}^n)$
under $\rho$ as a very proximal transformation, whose repelling
neighborhoods lie inside the $\epsilon$-repelling neighborhood of $g$ and
whose attracting points lies inside the $\epsilon$-attracting neighborhood
of $g$. We can certainly assume that $\rho (\delta_{m_1})$ satisfies the
conditions of Lemma \ref{fix}. Hence $\delta_{m_1}$ fixes some attracting
points $\overline{v}_{\delta_{m_1}},\overline{v}_{\delta_{m_1}^{-1}}$ which
are close to, but distinct from $\overline{v}_g,\overline{v}_{g^{-1}}$
respectively. Similarly the repelling neighborhoods of $\delta_{m_1},%
\delta_{m_1}^{-1}$ lie inside the $\epsilon$-repelling neighborhood of $%
g,g^{-1}$, and the repelling hyperplanes $\overline{H}_{\delta_{m_1}},%
\overline{H}_{\delta_{m_1}^{-1}}$ are close to that of $g$. We claim that
for all large enough $m_1$ 
\begin{equation}  \label{999}
\{\overline{v}_g,\overline{v}_{g^{-1}}\}\cap\big(\overline{H}%
_{\delta_{m_1}}\cup\overline{H}_{\delta_{m_1}^{-1}}\big)=\emptyset, \text{%
~and~} \{\overline{v}_{\delta_{m_1}},\overline{v}_{\delta_{m_1}^{-1}}\}\cap%
\big(\overline{H}_g\cup\overline{H}_{g^{-1}}\big)=\emptyset.
\end{equation}

Let us explain, for example, why $\overline{v}_{g^{-1}}\notin\overline{H}%
_{\delta_{m_1}}$ and why $\overline{v}_{\delta_{m_1}}\notin\overline{H}%
_{g^{-1}}$ (the other six conditions are similarly verified). Apply $%
\delta_{m_1}$ to the point $\overline{v}_{g^{-1}}$. As $g$ stabilizes $%
\overline{v}_{g^{-1}}$ we see that 
\begin{equation*}
\delta_{m_1}(\overline{v}_{g^{-1}})=g^{m_1}\gamma g^{m_1}(\overline{v}%
_{g^{-1}})=g^{m_1}\gamma (\overline{v}_{g^{-1}}).
\end{equation*}
Now, by our assumption, $\gamma (\overline{v}_{g^{-1}})\notin\overline{H}_g$%
. Moreover when ${m_1}$ is large, $g^{m_1}$ is a $\epsilon_{m_1}$%
-contracting with $\overline{H}_{g^{m_1}}=\overline{H}_g,~\overline{v}%
_{g^{m_1}}=\overline{v}_g$ and $\epsilon_{m_1}$ arbitrarily small. Hence, we
may assume that $\gamma (\overline{v}_{g^{-1}})$ is outside the $%
\epsilon_{m_1}$ repelling neighborhood of $g^{m_1}$. Hence $\delta_{m_1}(%
\overline{v}_{g^{-1}})=g^{m_1}\big(\gamma (\overline{v}_{g^{-1}})\big)$ lie
near $\overline{v}_g$ which is far from $\overline{H}_{\delta_{m_1}}$. Since 
$\overline{H}_{\delta_{m_1}}$ is invariant under $\delta_{m_1}$, we conclude that 
$\delta_{m_1}\overline{v}_{g^{-1}}\notin\overline{H}_{\delta_{m_1}}$.

To show that $\overline{v}_{\delta_{m_1}}\notin\overline{H}_{g^{-1}}$ we
shall apply $g^{-2{m_1}}$ to $\overline{v}_{\delta_{m_1}}$. If ${m_1}$ is
very large then $\overline{v}_{\delta_{m_1}}$ is very close to $\overline{v}%
_g$, and hence also $g^{m_1}(\overline{v}_{\delta_{m_1}})$ is very close to $%
\overline{v}_g$. As we assume that $\gamma$ takes $\overline{v}_g$ outside $%
\overline{H}_{g^{-1}}$, we get (by taking ${m_1}$ sufficiently large) that $%
\gamma$ also takes $g^{m_1}\overline{v}_{\delta_{m_1}}$ outside $\overline{H}%
_{g^{-1}}$. Taking ${m_1}$ even larger if necessary we get that $g^{-{m_1}}$
takes $\gamma g^{m_1}\overline{v}_{\delta_{m_1}}$ to a small neighborhood of 
$\overline{v}_{g^{-1}}$. Hence 
\begin{equation*}
g^{-2{m_1}}\overline{v}_{\delta_{m_1}}=g^{-2{m_1}}\delta_{m_1}\overline{v}%
_{\delta_{m_1}}=g^{-{m_1}}\gamma g^{m_1}\overline{v}_{\delta_{m_1}}
\end{equation*}
lies near $\overline{v}_{g^{-1}}$. Since $\overline{H}_{g^{-1}}$ is $g^{-2{%
m_1}}$ invariant and is far from $\overline{v}_{g^{-1}}$, we conclude that $%
\overline{v}_{\delta_{m_1}}\notin\overline{H}_{g^{-1}}$.

Now it follows from (\ref{999}) and Lemma \ref{fix} that for every $%
\epsilon_1>0$ we can take $j_1$ sufficiently large so that $g^{j_1}$ and $%
\delta_{m_1}^{j_1}$ are $\epsilon_1$-very proximal transformations, and the $%
\epsilon_1$-repelling neighborhoods of each of them are disjoint from the $%
\epsilon_1$-attracting points of the other, and hence they form a ping-pong
pair. Set $g_1=\delta_{m_1}^{j_1}$.

In a second step, we construct $g_2$ in an analogous way to the first step,
working with $g^{j_1}$ instead of $g$. In this way we would get $g_2$ which
is $\epsilon_2$-very proximal, and play ping-pong with $g^{j_1j_2}$.
Moreover, by construction, the $\epsilon_2$-repelling neighborhoods of $g_2$
lie inside the $\epsilon_1$-repelling neighborhoods of $g^{j_1}$, and the $%
\epsilon_2$-attracting neighborhoods of $g_2$ lie inside the $\epsilon_1$%
-attracting neighborhoods of $g^{j_1}$. Hence the three elements $g_1,~g_2$
and $g^{j_1j_2}$ form a ping-pong 3-tuple.

We continue recursively and construct the desired sequence $(g_n)$. Note
that in the Archimedean case, we have to make sure that the $g_n$'s form a
ping-pong $\aleph_0$-tuple also when we add to them the finitely many $z_i$%
's. This can be done by declaring $g_i=z_i$ for $i=1,\ldots l$, and starting
the recursive argument by constructing $g_{l+1}$.

Now since all $\Gamma _{k_{j}}=G_{k_{j}}\cap \Gamma $'s$\,$are mapped under
the homomorphism $\pi $ to Zariski dense subsets of $\Bbb{H}^{\circ }$, we
can multiply $x_{j}$ on the left and on the right by some elements of $%
\Gamma _{k_{j}}$ so that, if we call this new element $x_{j}$ again, $\rho
(x_{j})\overline{v}_{g_{j}}\notin \overline{H}_{g_{j}}$ and $\rho
(x_{j}^{-1})\overline{v}_{g_{j}^{-1}}\notin \overline{H}_{g_{j}^{-1}}$.
Considering the element $y_{j}=g_{j}^{l_{j}}x_{j}g_{j}^{l_{j}}$ for some
positive power $l_{j}$, we see that it lies in $\Gamma _{k_{j}}x_{j}\Gamma
_{k_{j}}$. Moreover, if we take $l_{j}$ large enough, it will behave on $%
\Bbb{P}(V_{K})$ like a very proximal transformation whose attracting and
repelling neighborhoods are contained in those of $g_{j}$. Therefore, the $%
y_{j} $'s also form an infinite ping-pong tuple and in the Archimedean case
they do so together with $z_{1},\ldots,z_{l}$. Hence the family $(y_{j})_{j}$
(resp. $(z_{1},\ldots,z_{l},(y_{j})_{j})$) generates a free group.

In the non-Archimedean case, the $y_{j}$'s are already dense in $G$ since we
selected them from sets which form a base for the topology. In the
Archimedean case, the elements $z_{1},\ldots ,z_{k}$ already generate a
dense subgroup of $H$, and since the $g_{j}$'s belong to $H$ and the $x_{j}$%
's are dense. Hence the group generated by the $z_{i}$'s and $y_{j}$'s
is dense in $G$. This completes the proof of Theorem \ref{infinite}.


\section{Applications to pro-finite groups}\label{applications}\label{profinite}

We derive two conclusions in the theory of pro-finite groups. The following
was conjectured by Dixon, Pyber, Seress and Shalev (see \cite{DPSS}):

\begin{thm}
\label{dpss}Let $\Gamma $ be a finitely generated linear group over some
field. Assume that $\Gamma $ is not virtually solvable. Then, for any
integer $r\geq d(\Gamma )$, its pro-finite completion $\hat{\Gamma}$
contains a dense free subgroup of rank $r.$
\end{thm}

\proof%
Let $R$ be the ring generated by the matrix entries of the elements of $%
\Gamma $. It follows from the Noether normalization theorem that $R$ can be
embedded in the valuation ring ${\mathcal{O}}$ of some local field $k$. Such
an embedding induces an embedding of $\Gamma $ in the pro-finite group $%
\text{GL}_{n}({\mathcal{O}})$. By the universal property of $\hat{\Gamma}$
this embedding induces a surjective map $\hat{\Gamma}\to \overline{\Gamma }%
\leq \text{GL}_{n}({\mathcal{O}})$ onto the closure of the image of $\Gamma $
in $\text{GL}_{n}(k)$. Since $\Gamma $ is not virtually solvable, $\overline{%
\Gamma }$ contains no open solvable subgroup, and hence by Theorem \ref{main}
(see also Proposition \ref{gC-compact}), $\overline{\Gamma }$ contains a
dense $F_{r}$ (in fact we can find such an $F_{r}$ inside $\Gamma $). By
Gasch\"{u}tz's lemma (see \cite{Rib}, Proposition 2.5.4) it is possible to
lift the $r$ generators of this $F_{r}$ to $r$ elements in $\hat{\Gamma}$
generating a dense subgroup in $\hat{\Gamma}$. These lifts, thus, generate a
dense $F_{r}$ in $\hat{\Gamma}$. 
\endproof%

\medskip

Let now $H$ be a subgroup of a group $G$. Following \cite{Sh} we define the
notion of coset identity as follows:

\begin{defn}
A group $G$ satisfies a coset identity with respect to $H$ if there exist

\begin{itemize}
\item  a non-trivial reduced word $W$ on $l$ letters,

\item  $l$ fixed elements $g_{1},\ldots ,g_{l}$,
\end{itemize}

such that the identity 
\begin{equation*}
W(g_{1}h_{1},\ldots ,g_{l}h_{l})=1
\end{equation*}
holds for any $h_{1},\ldots ,h_{l}\in H$.
\end{defn}

It was conjectured by Shalev \cite{Sh} that if there is a coset identity
with respect to some open subgroup in a pro-$p$ group $G$, then there is
also an identity in $G$. The following immediate consequence of Corollary 
\ref{freecosets} settles this conjecture in the case where $G$ is an
analytic pro-$p$ group, and in fact, shows that a stronger statement is true
in this case:

\begin{thm}
Let $G$ be an analytic pro-$p$ group. If $G$ satisfies a coset identity with
respect to some open subgroup, then $G$ is virtually solvable.
\end{thm}

\proof
If $G$ is not virtually solvable, then by Corollary \ref{freecosets}, we can
choose coset representatives for $H$ in $G$ which are free generators of a
free group. \endproof

In fact the analogous statement holds also for finitely generated linear
groups:

\begin{thm}
Let $\Gamma $ be a finitely generated linear group over any field. If $%
\Gamma $ satisfies a coset identity with respect to some finite index
subgroup $\Delta $, then $\Gamma $ is virtually solvable.
\end{thm}


\section{Applications to amenable actions}
\label{amenable}

For convenience, we introduce the following definition:

\begin{defn}
We shall say that a topological group $G$ has property (OS) if it contains
an open solvable subgroup.
\end{defn}

Our main result, Theorem \ref{main}, states that if $\Gamma $ is a (finitely
generated) linear topological group over a local field, then either $\Gamma $
has property (OS) or $\Gamma $ contains a dense (finitely generated) free
subgroup. In the previous section we proved the analogous statement for
pro-finite completions of linear groups over an arbitrary field. For real
Lie groups, property (OS) is equivalent to ``the identity component is
solvable''.

It was conjectured by Connes and Sullivan and proved subsequently by Zimmer 
\cite{Zim} that if $\Gamma $ is a countable subgroup of a real Lie group $G$%
, then the action of $\Gamma $ on $G$ by left multiplications is amenable if
and only if $\Gamma $ has property (OS). Note also that if $\gC$ acts amenably on $G$, then it also acts amenably on $G/P$
whenever $P\leq G$ is closed amenable subgroup.
We refer the reader to \cite{Zim2}
Chapter 4 for an introduction and background on amenable actions. The harder
part of the equivalence is to show that if $\Gamma $ acts amenably then it
has (OS). As noted by Carri\`{e}re and Ghys \cite{Ghys}, the Connes-Sullivan
conjecture is a straightforward consequence of Theorem \ref{main}. Let us
reexplain this claim: by Theorem \ref{main}, it is enough to show that if $%
\Gamma $ contains a non-discrete free subgroup, then it cannot act amenably.

\proof%
(\textit{non-discrete free subgroup}$\Rightarrow $\textit{action is
non-amenable}). By contradiction, if $\Gamma $ were acting amenably, then
any subgroup would do so too, hence we can assume that $\Gamma $ itself is a
non-discrete free group $\left\langle x,y\right\rangle .$ By Proposition
4.3.9 in \cite{Zim2}, it follows that there exists a $\Gamma $-equivariant
Borel map $g\mapsto m_{g}$ from $G$ to the space of probability measures on
the boundary $\partial \Gamma $. Let $X$ (resp. $Y$) be the set of infinite
words starting with a non trivial power of $x$ (resp. $y$). Let $(\xi _{n})$
(resp. $\theta _{n}$) be a sequence of elements of $\Gamma $ tending to the
identity element in $G$ and consisting of reduced words starting with $y$
(resp. $y^{-1}$). By the converse to Lebesgue's dominated convergence
theorem, up to passing to a subsequence of $(\xi _{n})_{n}$ if necessary, we
have that for almost all $g\in G$, $m_{g\xi _{n}}(X)$ and $m_{g\theta
_{n}}(X)$ converge to $m_{g}(X)$. However, for almost every $g\in G$, $%
m_{g\xi _{n}}(X)=m_{g}(\xi _{n}X)$ and $m_{g\theta _{n}}(X)=m_{g}(\theta
_{n}X)$. Moreover $\xi _{n}X$ and $\theta _{n}X$ are disjoint subsets of $Y$%
. Hence, for almost every $g\in G$, $2m_{g}(X)\leq m_{g}(Y)$. Reversing the
roles of $X$ and $Y$ we get a contradiction. 
\endproof%

Clearly, this proof is valid whenever $G$ is a locally compact group. In particular, Theorem \ref{main}
implies the following general result:

\begin{thm}
Let $k$ be a local field, $G$ a closed subgroup of $\GL_n(k)$, $P$ a closed amenable subgroup of $G$ and 
$\gC\leq G$ a countable subgroup. Then the following are equivalent:
\begin{enumerate}
\item
$\gC$ has property (OS).
\item
The action of $\gC$ by left multiplications on the homogeneous space $G/P$ is amenable.
\item
$\gC$ contains no non-discrete free subgroup.
\end{enumerate}
\end{thm}

A theorem of Auslander (see \cite{raghunathan} 8.24) states that if $G$ is a
real Lie group, $R$ a closed normal solvable subgroup, and $\Gamma $ a
subgroup with property (OS), then the image of $\Gamma $ in $G/R$ also has
property (OS). Taking $G$ to be the group of Euclidean motions, $R$ the
subgroup of translations, and $\Gamma \leq G$ a torsion free lattice, one
obtains the classical theorem of Bieberbach that any compact Euclidean
manifold is finitely covered by a torus.

Following Zimmer (\cite{Zim}), we remark that Auslander's theorem follows from Zimmer's theorem. To see
this, note that $G$ (hence also $\Gamma $) being second countable, we can
always replace $\Gamma $ by a countable dense subgroup of it. Then, if $%
\Gamma $ has property (OS) it must act amenably on $G$. As $R$ is closed and
amenable, this implies that $\Gamma R/R$ acts amenably on $G/R$ (see \cite
{Zim2} chapter 4), which in turn implies, by Zimmer's theorem, that $\Gamma R/R$ has
property (OS).

A discrete linear group is amenable if and only if it is virtually solvable.
It follows that for a countable linear group over some topological field,
being (OS) is the same as containing an open amenable subgroup.

\begin{defn}
We shall say that a countable topological group $\Gamma $ has property (OA)
if it contains an open subgroup, which is amenable in the abstract sense
(i.e. amenable with respect to the discrete topology).
\end{defn}

The following is a generalization of Zimmer's theorem:

\begin{thm}
\label{G-CS} Let $G$ be a locally compact group, and let $\Gamma \leq G$ be
a countable subgroup. Then the action of $\Gamma $ on $G$ by left
multiplications is amenable if and only if $\Gamma $ has property (OA).
\end{thm}

\begin{proof} 
The proof makes use of the structure theory
for locally compact groups (see \cite{MZ}). We shall reduce the general case
to the already known case of real Lie groups.

The ``if'' side is clear.

Assume that $\Gamma $ acts amenably. Let $G^{0}$ be the identity connected
component of $G$. Then $G^{0}$ is normal in $G$ and $F=G/G^{0}$ is a totally
disconnected locally compact group, and as such, has an open profinite
subgroup. Since $\Gamma $ acts amenably, its intersection with an open
subgroup acts amenably on the open subgroup. Therefore we can assume that $%
G/G^{0}$ itself is profinite. By \cite{MZ} Theorem 4.6, there is a compact
normal subgroup $K$ in $G$ such that the quotient $G/K$ is a Lie group. Up
to passing to an open subgroup of $G$ again, we can assume that $G/K$ is
connected. Since $\Gamma \cap K$ acts amenably on $K$ and $K$ is amenable, $%
\Gamma \cap K$ is amenable (see \cite{Zim2}, Chapter 4). Moreover, as $K$ is
amenable, $\Gamma $, and hence also $\Gamma K/K$, acts amenably on the
connected Lie group $G/K$. We conclude that $\Gamma K/K$ has property (OS)
and hence $\Gamma $ has property (OA). 
\end{proof}

As an immediate corollary we obtain the following generalization of
Auslander's theorem:

\begin{thm}
\label{G-Aus} Let $G$ be locally compact group, $R\leq G$ a closed normal
amenable subgroup, and $\Gamma \leq G$ a subgroup with property (OA). Then
the image of $\Gamma $ in $G/R$ has also (OA).
\end{thm}

\begin{rem}
The original statement of Auslander follows easily from \ref{G-Aus}.
\end{rem}

As a consequence of Theorem \ref{G-Aus} we derive a structural result for lattices in general locally 
compact groups. Let $G$ be a locally compact group. Then $G$ admits a unique maximal closed normal amenable 
subgroup $P$, and $G/P$ is isomorphic (up to finite index) to a direct product 
$$
G\cong G_d\times G_c,
$$
of a totally disconnected group $G_d$ with a connected center-free semisimple Lie group without compact 
factors $G_c$ (see \cite{BurgerMonod} Theorem 3.3.3). Let $\gC\leq G$ be a lattice, then we have:

\begin{prop}[This result was proved in a conversation with Marc Burger]\label{G_cxG_d}
The projection of $\gC$ to the connected factor $G_c$ lies between a lattice in $G_c$ to its commensurator.
\end{prop}

\begin{proof}
Let $\pi,\pi_d,\pi_c$ denote the quotient maps from $G$ to $G/P,G_d,G_c$ respectively.
Since $\gC$ is discrete, and hence has property (OA) it follows from Theorem \ref{G-Aus} that 
$\gD=\pi (\gC )$ also has (OA). Let $A\leq G_d$ be an open compact subgroup, and let 
$\gD^0=\gD\cap (A\times G_c)$. Then $\gD^0$ has (OA) being open in $\gD$.
Clearly, $\gD$ commensurates $\gD^0$.
Let $\gS$ be the projection of $\gD^0$ to $G_c$. By Theorem \ref{G-Aus} $\gS$ has (OA), and hence also 
$\overline{\gS}$ has (OA), i.e. the identity connected component $\overline{\gS}^0$ is solvable.
However, since the homogeneous space $G_c/\overline{\gS}$ clearly carries a finite $G_c$-invariant Borel 
regular measure, it follows from Borel's density theorem that  
$\overline{\gS}^0$ is normal in $G_c$. 
Since $\overline{\gS}^0$ is also solvable and $G_c$ is semisimple, it follows that $\overline{\gS}$
is discrete. Therefore $\overline\gS =\gS$ and $\gS$ is a lattice in $G_c$, and 
$\gS\leq\pi_c(\gC)\leq Comm_{G_c}(\gS )$.
\end{proof}


\section{The growth of leaves in Riemannian foliations}\label{growth}

The main result of this section is the following theorem which answers a
question of Carri\`{e}re \cite{Car1} (see also \cite{Hae1}):

\begin{thm}
\label{youyou}Let $\mathcal{F}$ be a Riemannian foliation on a compact
manifold $M$. The leaves of $\mathcal{F}$ have polynomial growth if and only
if the structural Lie algebra of $\mathcal{F}$ is nilpotent. Otherwise, generic leaves have exponential growth.
\end{thm}

In fact, in the second case, we show that the holonomy cover of an arbitrary leaf has exponential volume growth. The result then follows from the fact that a generic leaf has no holonomy (see \cite{Con}). The example below show that it is possible that some leaves are compact while generic leaves have exponential growth.
For background and definitions about Riemannian foliations, the structural
Lie algebra, and growth of leaves see \cite{Con} \cite{Car1}, \cite{Hae1} and \cite{Mol}

\begin{exam}
We give here an example of a Riemannian foliation on a compact manifold which has two compact leaves although every other leaf has exponential volume growth. This example was shown to us by E. Ghys. Let $\Gamma$ be the fundamental group of a surface of genus $g\geq2$ and $\pi : \Gamma \rightarrow SO(3)$ be a faithful representation. Such a map can be obtained, for instance, by realizing $\Gamma$ as a torsion free co-compact arithmetic lattice in $SO(2,1)$. The group $\Gamma$ acts freely and co-compactly on the hyperbolic plane $\Bbb{H}^2$ by deck transformations and it acts via the homomorphism $\pi$ by rotations on the $3$-sphere $S^3$ viewed as $\Bbb{R}^3\cup\{\infty\}$. We now \textit{suspend} $\pi$ to the quotient manifold $M=(\Bbb{H}^2 \times S^3)/\Gamma$ and obtain a Riemannian foliation on $M$ whose leaves are projections to $M$ of each $(\Bbb{H}^2,y)$, $y \in S^3$. There are three types of leaves. The group $SO(3)$ has two fixed points on $S^3$, the north and south poles, each giving rise to a compact leaf $\Bbb{H}^2/\Gamma$ in $M$. Each $\gamma \in \Gamma \setminus \{0\}$ stabilizes a line in $S^3$ and the leaf through each point of this line other than the poles will be a hyperbolic cylinder $\Bbb{H}^2/\langle \gamma \rangle$. Any other point in $S^3$ has a trivial stabilizer in $\Gamma$ and gives rise to a leaf isometric to $\Bbb{H}^2$. These are the generic holonomy-free leaves. Apart from the two compact leaves all others have exponential volume growth.   
\end{exam}

Following Carri\`{e}re \cite{Car1}, \cite{DCa} we define the \textbf{local
growth} of a finitely generated subgroup $\Gamma $ in a given connected real
Lie group $G$ in the following way. Fix a left-invariant Riemannian metric
on $G$ and consider the open ball $B_{R}$ of radius $R>0$ around the
identity. Suppose that $S$ is a finite symmetric set of generators of $%
\Gamma .$ Let $B(n)$ be the ball of radius $n$ in $\Gamma $ for the word
metric determined by $S$, and let $B_{R}(n)$ be the subset of $B(n)$
consisting of those elements $\gamma \in B(n)$ which can be written as a
product $\gamma =\gamma _{1}\cdot \ldots \cdot \gamma _{k}$, $k\leq n$, of
generators $\gamma _{i}\in S$ in such a way that whenever $1\leq i\leq k$
the element $\gamma _{1}\cdot \ldots \cdot \gamma _{i}$ belongs to $B_{R}$.
In this situation, we say that $\gamma $ can be written as a word with
letters in $S$ which \textit{stays all its life in }$B_{R}$. Let $%
f_{R,S}(n)=card(B_{R}(n))$. As it is easy to check, if $S_{1}$ and $S_{2}$
are two symmetric sets of generators of $\Gamma $, then there exist integers 
$N_{0},N_{1}>0$ such that $f_{R,S_{1}}(n)\leq f_{R+N_{0},S_{2}}(N_{1}n)$.

\begin{defn}\label{local-growth}
The local growth of $\Gamma $ in $G$ with respect to a set $S$ of generators
and a ball $B_{R}$ of radius $R$ is the growth type of $f_{R,S}(n)$.
\end{defn}

\noindent The growth type of $f_{R,S}(n)$ is \textit{polynomial} if there
are positive constants $A$ and $B$ such that $f_{R,S}(n)\leq An^{B}$ and is 
\textit{exponential} if there are constants $C>0$ and $\rho >1$ such that $%
f_{R,S}(n)\geq C\rho ^{n}$. It can be seen that $\Gamma $ is discrete in $G$
if and only if the local growth is bounded for any $S$ and $R$.

According to Molino's theory (see \cite{Hae1}, \cite{Mol}) every Riemannian foliation $\mathcal{F}$ on a compact manifold $M$ lifts to a foliation $\mathcal{F}'$ of the same dimension on the bundle $M'$ of transverse orthonormal frames over $M$. Moreover, when restricting the foliation $\mathcal{F'}$ to the closure of a given leaf, we obtain a Lie foliation for some simply connected Lie group $G$, the structural Lie group. This allows to reduce the problem to Lie foliations with dense leaves. A $G$-Lie foliation is by definition a foliation where local transverse manifolds are identified with open subsets of $G$ and transitions maps are given by left translations in $G$. Given a base point $x_0$ in $M$ we obtain a natural homomorphism of the fundamental group of $M$ into $G$, whose image is called the holonomy group of the Lie foliation and is a dense subgroup $\Gamma$ in $G$. As Carri\`ere pointed out in \cite{Car1}, the volume growth of any leaf of a $G$-Lie foliation is coarsely equivalent to the local growth of $\Gamma$ in $G$. Hence Theorem \ref{youyou} is a
consequence of the following:

\begin{thm}
\label{local}Let $\Gamma $ be a finitely generated dense subgroup of a
connected real Lie group $G$. If $G$ is nilpotent then $\Gamma $ has
polynomial local growth (for any choice of $S$ and $R$). If $G$ is not
nilpotent, then $\Gamma $ has exponential local growth (for any choice of $S$
and any $R$ big enough).
\end{thm}

The rest of this section is therefore devoted to the proof of Theorem \ref
{local}. As it turns out, Theorem \ref{local} is an immediate corollary of
Theorem \ref{main} in the case when $G$ is not solvable. When $G$ is
solvable we can adapt the argument as shown below. The main proposition is
the following:

\begin{prop}
\label{perturb}\label{freesemi} Let $G$ be a non-nilpotent connected real
Lie group and $\Gamma $ a finitely generated dense subgroup. For any finite
set $S=\{s_{1},\ldots ,s_{k}\}$ of generators of $\Gamma $, and any $%
\varepsilon >0$, one can find perturbations $t_{i}\in \Gamma $ of the $%
s_{i},~i=1,\ldots ,k$ such that $t_{i}\in s_{i}B_{\varepsilon }$ and the $%
t_{i}$'s are free generators of a free semi-group on $k$ generators.
\end{prop}

Before going through the proof of this proposition, let us explain how we
deduce from it a proof of Theorem \ref{local}.

\begin{proof}[Proof of Theorem \ref{local}.] Suppose that $\Sigma :=\{g_{1},\ldots
,g_{N},h_{1},\ldots ,h_{N}\}$ is a subset of $B_{R}$ consisting of pairwise
distinct elements such that both $\{g_{1},\ldots ,g_{N}\}$ and $%
\{h_{1},\ldots ,h_{N}\}$ are maximal $R/2$-discrete subsets of $\overline{%
B_{R}}$ (that is $d(g_{i},g_{j}),d(h_{i},h_{j})\geq R/2$ if $i\neq j$). Then 
\begin{equation*}
\overline{B_{R}}\subset \bigcup_{1\leq i,j\leq N}(g_{i}B_{R/2}\cap
h_{j}B_{R/2}).
\end{equation*}

\begin{lem}\label{ammel} 
Let $G$ be a connected real Lie group endowed with a
left-invariant Riemannian metric. Let $B_{R}$ be the open ball of radius $R$ centered
at the identity. Let $\Sigma =\{s_{1},\ldots ,s_{k}\}$ be a finite subset of
pairwise distinct elements of $B_{R}$ such that 
\begin{equation}
\overline{B_{R}}\subset \bigcup_{i<j}(s_{i}^{-1}B_{R/2}\cap
s_{j}^{-1}B_{R/2}).  \label{cover}
\end{equation}
Assume also that the elements of $\Sigma $ are free generators of a free
semi-group. Then any finitely generated subgroup of $G$ containing $\Sigma $
has exponential local growth.
\end{lem}

\begin{proof}
Let $S(n)$ be the sphere of radius $n$ in the free
semi-group for the word metric determined by the generating set $\Sigma $.
Let $w\in S(n)\cap B_{R}$. By $(\ref{cover})$ there are
indices $i\neq j$ such that $w\in s_{i}^{-1}B_{R/2}$ and $w\in
s_{j}^{-1}B_{R/2}$. This implies that $s_{i}w$ and $s_{j}w$ belong to $%
S(n+1)\cap B_{R}(n+1)$. All elements obtained in this way are pairwise distinct,
hence $card(S(n+1)\cap B_{R}(n+1))\geq 2\cdot card(S(n)\cap B_{R}(n))$. This yields $%
card(S(n)\cap B_{R}(n))\geq 2^{n}$ for all $n\geq 0$. 
\end{proof}

Now observe that any small enough perturbation of a finite set $\Sigma $ in $%
B_{R}$ satisfying $(\ref{cover})$ still satisfies $(\ref{cover}).$ Hence
exponential local growth for dense subgroups in non-nilpotent connected real
Lie groups follows from the combination of Lemma \ref{ammel} and Proposition 
\ref{perturb}. 
\end{proof}

\textit{Proof of Proposition \ref{perturb}.} When $G$ is not solvable, we
already know this fact from the proof of Theorem \ref{main} for connected
Lie groups (see Paragraph \ref{archi}). In that case, we showed that we
could even take the $t_{i}$'s to generate a free subgroup. Thus we may
assume that $G$ is solvable. By Ado's theorem it is locally isomorphic to a
subgroup of $\text{GL}_{n}({\mathbb{C}})$, and it is easy to check that the
property to be shown in Proposition \ref{perturb} does not change by local
isomorphisms. Thus, we may also assume that $G\leq\text{GL}_{n}(\Bbb{C})$.
Let $\Bbb{G}$ be the Zariski closure of $G$ in $GL_{n}(\Bbb{C})$. It is a
Zariski connected solvable algebraic group over $\Bbb{C}$ which is not
nilpotent. We need the following elementary lemma for $k=\Bbb{C}$.

\begin{lem}
\label{affine}Let $\Bbb{G}$ be a solvable connected algebraic $k$-group
which is not nilpotent. Suppose it is $k$-split, then there is an algebraic $%
k$-morphism from $\Bbb{G(}k\Bbb{)}$ to $\Bbb{GL}_{2}(k)$ whose image is the
full affine group 
\begin{equation}
\Bbb{A}(k)=\left\{ \left( 
\begin{array}{ll}
a & b \\ 
0 & 1
\end{array}
\right) \,a,b\in k\right\} .  \label{affineg}
\end{equation}
\end{lem}

\begin{proof}
We proceed by induction on $\dim \Bbb{G}$. We can write $G:=\Bbb{G}%
(k)=T\cdot N$ where $T=\Bbb{T(}k\Bbb{)}$ is a split torus and $N=\Bbb{N(}k%
\Bbb{)}$ is the unipotent radical of $G$ (see \cite{Bo}, Chapter III)$.$ Let 
$Z$ be the center of $N$. It is a non trivial normal algebraic subgroup of $G
$. If $T$ acts trivially on $Z$ by conjugation then $G/Z$ is again
non-nilpotent $k$-split solvable $k$-group and we can use induction. thus we
may assume that $T$ acts non-trivially on $Z$ by conjugation. As $T$ is
split over $k$, its action on $Z$ also splits, and there is a non-trivial
algebraic multiplicative character $\chi :T\to \Bbb{G}_{m}(k)$ defined over $%
k$ and a $1$-dimensional subgroup $Z_{\chi }$ of $Z$ such that, identifying $%
Z_{\chi }$ with the additive group $\Bbb{G}_{a}(k)$, we have $tzt^{-1}=\chi
(t)z$ for all $t\in T$ and $z\in Z_{\chi }$. It follows that $Z_{\chi }$ is
a normal subgroup in $G$, and we can assume that $T$ acts trivially on $%
N/Z_{\chi }$, for otherwise we could apply the induction assumption on $%
G/Z_{\chi }$. For all $\gamma \in T$, this yields a homomorphism $\pi
_{\gamma }:N\rightarrow Z_{\chi }$ given by the formula $\pi _{\gamma
}(n)=\gamma n\gamma ^{-1}n^{-1}$. Since $T$ and $N$ do not commute, $\pi
_{\gamma }$ is non trivial for at least one $\gamma \in T$. Fix such a $%
\gamma $ and let $N$ act on $Z_{\chi }$ by left multiplication by $\pi
_{\gamma }(n)$. Let $T$ act on $Z_{\chi }$ by conjugation. One can verify
that this yields an algebraic action of the whole of $G$ on $Z_{\chi }$.
Identifying $Z_{\chi }$ with the additive group $\Bbb{G}_{a}(k),$ we have
that $N$ acts unipotently and non-trivially and $T$ acts via the non-trivial
character $\chi $. We have found a $k$-algebraic affine action of $G$ on the
line, and hence a $k$-map ${\mathbb{G}}\to {\mathbb{A}}$. Clearly this map
is onto. {$\square $ }

By Lemma \ref{affine}, $\Bbb{G}$ surjects onto the affine group of the
complex line, which we denote by $A=\Bbb{A}(\Bbb{C})$. The image of $G$ is a
real connected subgroup of $A$ which is Zariski dense. Hence it is enough to
prove Proposition \ref{freesemi} for Zariski dense connected subgroups of $A$%
. We need the following technical lemma:

\begin{lem}
Let $\Gamma $ be a non-discrete finitely generated Zariski dense subgroup of 
$\Bbb{A}({\mathbb{C}})$ with connected closure. Let $R\subset {\mathbb{C}}$
be the subring generated by the matrix entries of elements in $\Gamma $.
Then there exists a sequence $(\gamma _{n})_{n}$ of points of $\Gamma ,$
together with a ring embedding $\sigma :R\hookrightarrow k$ into another
local field $k$, such that $\gamma _{n}=(a_{n},b_{n})\rightarrow (1,0)$ in $%
\Bbb{A}(\Bbb{C})$ and $\sigma (\gamma _{n})=\big(\sigma (a_{n}),\sigma
(b_{n})\big)\rightarrow (0,\sigma (\beta ))$ in the topology of $k$ for some
number $\beta $ in the field of fractions of $R$.
\end{lem}

\begin{proof}

Let $g_{n}=(a_{n},b_{n})$ be a sequence of distinct elements of $\Gamma $
converging to identity in $\Bbb{A}({\mathbb{C}})$ and such that $|a_{n}|_{{%
\mathbb{C}}}\leq 1$ and $a_{n}\neq 1$ for all integers $n$. From Lemma \ref
{Tit} one can find a ring embedding $\sigma :R\hookrightarrow k$ for some
local field $k$ such that, up to passing to a subsequence of $g_{n}$'s, we
have $\sigma (a_{n})\rightarrow 0$ in $k$. We can assume $|\sigma
(a_{n})|_{k}<1$ for all $n$. Now let $\xi =(a,b):=g_{0}$ and consider the
element 
\begin{equation*}
\xi ^{m}g_{n}\xi ^{-m}=(a_{n},\frac{1-a^{m}}{1-a}b(1-a_{n})+a^{m}b_{n}).
\end{equation*}
Since $|a|_{{\mathbb{C}}}\leq 1$, the second component remains $\leq \frac{2%
}{|1-a|_{{\mathbb{C}}}}|b|_{{\mathbb{C}}}|1-a_{n}|_{{\mathbb{C}}}+|b_{n}|_{{%
\mathbb{C}}}$ for all $m$, and tends to $0$ in ${\mathbb{C}}$ when $n\to
\infty $ uniformly in $m$. Applying the isomorphism $\sigma $, we have: 
\begin{equation}
\sigma (\xi ^{m}g_{n}\xi ^{-m})=\big(\sigma (a_{n}),\frac{1-\sigma (a)^{m}}{%
1-\sigma (a)}\sigma (b)(1-\sigma (a_{n}))+\sigma (a)^{m}\sigma (b_{n})\big).
\label{eq22}
\end{equation}
Since $|\sigma (a)|_{k}<1$, for any given $n$, choosing $m$ large, we can
make $|\sigma (a)^{m}\sigma (b_{n})|_{k}$ arbitrarily small. Hence for some
sequence $m_{n}\rightarrow +\infty $ the second component in (\ref{eq22})
tends to $\sigma (\beta )$ where $\beta :=\frac{b}{1-a}$ as $n$ tends to $%
+\infty $. \end{proof}

We shall now complete the proof of Proposition \ref{perturb}. Note that if $k
$ is some local field and $\gamma =(a_{0},b_{0})\in \Bbb{A}(k)$ with $%
|a_{0}|_{k}<1$, then $\gamma$ acts on the affine line $k$ with a fixed 
point $%
x_{0}=b_{0}/(1-a_{0})$ and it contracts the disc of radius $R$ around $x_{0}$
to the disc of radius $|a_{0}|_{k}\cdot R$. Therefore, if we are given $t$
distinct points $b_{1},\ldots ,b_{t}$ in $k$, there exists $\varepsilon >0$
such that for all $a_{1},\ldots ,a_{t}\in k$ with $|a_{i}|_{k}\leq
\varepsilon $, $i=1,\ldots ,t$, the elements $(a_{i},b_{i})$'s play
ping-pong on the affine line, hence are free generators of a free
semi-group. The group $G=\overline{\Gamma }$ is a connected and Zariski
dense subgroup of $\Bbb{A}(\Bbb{C})$: it follows that we can find arbitrary
small perturbations $\widetilde{s_{i}}$ of the $s_{i}$'s within $\Gamma $
such that the $a(\tilde{s}_{i})\beta +b(\tilde{s}_{i})$'s are pairwise
distinct complex numbers. If $(\gamma _{n})_{n}$ is the sequence obtained in
the last Lemma, then for some $n$ large enough the points $t_{i}:=\widetilde{%
s_{i}}\gamma _{n}$ will be small perturbations of the $s_{i}$'s (i.e. belong
to $s_{i}B_{\varepsilon }$) and the $\sigma (t_{i})=\big(\sigma (a(\tilde{s}%
_{i})a_{n}),\sigma (a(\tilde{s}_{i})b_{n})+\sigma (b(\tilde{s}_{i}))\big)$
will play ping-pong on $k$ for the reason we just explained (the $\sigma (a(%
\tilde{s}_{i}))\sigma (\beta )+\sigma (b(\tilde{s}_{i}))$'s are all
distinct). \end{proof}

\begin{Ack}
We thank the following persons for their help and the interesting
conversations we had about various points in this article: Y. Barnea, M. Burger, D.
Getz, E. Ghys, M. Larsen, A. Lubotzky, G.A. Margulis, N. Nikolov, J.F.
Quint, A. Salehi-Golsefidi, G.A. Soifer. We also thank Olivia Barbarroux for
her hospitality during our stay at Luminy in July 2002 and H. Abels for
inviting us to Bielefeld in July 2003, where part of this work was conducted.
\end{Ack}

\end{document}